\documentclass[a4paper]{article}
\usepackage[english]{babel}
\usepackage[utf8]{inputenc}
\usepackage{amsmath,amsthm}
\usepackage{amsfonts,amssymb}

\usepackage{tikz,tikz-cd} 

\usepackage{xifthen} 

\newcommand{\toposdefaut}{0}
\newcommand{\topos}[1][\toposdefaut]{ 
\ifthenelse{\equal{#1}{0}}{ \mathcal{T} }
{
\ifthenelse{\equal{#1}{1}}{ \mathcal{E} }{ #1 }
}
}

\newcommand{\loc}{\text{Loc}}

\newcommand{\Z}{\mathbb{Z}}
\newcommand{\R}{\mathbb{R}}
\newcommand{\C}{\mathbb{C}}
\newcommand{\N}{\mathbb{N}}

\newcommand{\G}{\mathbb{G}}

\newcommand{\Ecal}{\mathcal{E}} 
 
\newcommand{\Tcal}{\mathcal{T}}

\newcommand{\Ocal}{\mathcal{O}}

\newcommand{\Scal}{\mathcal{S}} 
 
\newcommand{\Fcal}{\mathcal{F}}

\newcommand{\Kcal}{\mathcal{K}} 
\newcommand{\Lcal}{\mathcal{L}}

\newcommand{\Ccal}{\mathcal{C}} 
 
\newcommand{\Bcal}{\mathcal{B}}

\usepackage{titlesec}
\titleformat{\subsubsection}[runin]{\normalfont}{\thesubsubsection}{0pt}{}[.]

\renewcommand{\thesubsubsection}{\arabic{section}.\arabic{subsubsection}}
\csname @addtoreset\endcsname{subsubsection}{section}

\newcommand{\block}[1]
{

\par \subsubsection{} #1

\bigskip}

\newcommand{\blockn}[1]{\par #1 \bigskip}

\newcommand{\Th}[1]
	{
	\bigskip	
	\textbf{Theorem : }{\itshape #1}
		
	\bigskip
	}

\newcommand{\Prop}[1]
	{

	\bigskip
	
	\textbf{Proposition : }{\itshape #1}
		
	\bigskip
	
	}

\newcommand{\Cor}[1]
	{

	\bigskip
	
	\textbf{Corollary : }{\itshape #1}	
		
	\bigskip

	}

\newcommand{\Lem}[1]
	{

	\bigskip
	
	\textbf{Lemma : }{\itshape #1}
		
	\bigskip
	
	}

\newcommand{\Def}[1]
	{
	
	\bigskip
	
	\textbf{Definition : }{\itshape #1}
	
	\bigskip
	
	}

\newcommand{\Dem}[1]{
	
	\smallskip
	
	\textbf{Proof : } \par
	 {#1} $\square$
	 
	 \bigskip
}

\begin{document}

\pagestyle{plain}
\title{The convolution algebra of an absolutely locally compact topos}
\author{Simon Henry}

\maketitle

\begin{abstract}

We introduce a class of toposes called ``absolutely locally compact'' toposes and of ``admissible'' sheaf of rings over such toposes. To any such ringed topos $(\mathcal{T},A)$ we attach an involutive convolution algebra $\mathcal{C}_c(\mathcal{T},A)$ which is well defined up to Morita equivalence and characterized by the fact that the category of non-degenerate modules over $\mathcal{C}_c(\mathcal{T},A)$ is equivalent to the category of sheaf of $A$-module over $\mathcal{T}$. In the case where $A$ is the sheaf of real or complex Dedekind numbers, we construct several norms on this involutive algebra that allows to complete it in various Banach and $C^*$-algebras: $L^1(\Tcal,A)$, $C^*_{red}(\Tcal,A)$ and $C^*_{max}(\Tcal,A)$. We also give some examples where this construction corresponds to well known constructions of involutive algebras, like groupoids convolution algebra and  Leavitt path algebras.

\end{abstract}

\renewcommand{\thefootnote}{\fnsymbol{footnote}} 
\footnotetext{\emph{email:} simon.henry@college-de-france.fr}
\renewcommand{\thefootnote}{\arabic{footnote}} 


\tableofcontents

\section{Introduction}

\blockn{Both topos theory and non-commutative geometry are theories that are concerned with the study of certain ``generalized topological spaces'' that are too singular to be treated by ordinary topology, like the space of leaves of a foliation or the space of orbits of a group action or of a groupoid. Indeed $C^*$-algebras are supposed to describe some sort of generalized locally compact topological space while toposes are also generalized topological spaces. Moreover, a lot of tools available for ordinary topological spaces can be extended to one of theses generalized context, for example cohomology extend very naturally to toposes, measure theory and integration extend very naturally to $C^*$-algebras. There is also a lot of examples of objects to which one can attach both a $C^*$-algebra and a topos, for example foliations, dynamical systems and topological groupoids as mentioned above, but also graphs (and various generalization of graphs).}

\blockn{In this paper we will present a construction that attach to a topos satisfying certain conditions of local compactness and local separation an involutive convolution algebra (similar to the convolution algebra of compactly supported functions on a groupoid) that can then be completed into a reduced or a maximal $C^*$-algebra. In a large number of classical example of geometric object to which one can attach both a topos and a (reduced or maximal) $C^*$-algebra, the $C^*$-algebra can be recovered as the $C^*$-algebra attached to this topos by this construction.}

\blockn{More precisely, we start with a topos $\Tcal$ which is ``absolutely locally compact'' (see definition \ref{def_absloccpt}) endowed with an ``admissible'' (see definition \ref{DefAdmissibleSR} ) sheaf of rings $A$, then to any object of $\Tcal$ which is a bound and ``$A$-separating'' (see \ref{DefAdmissibleSR}) one can associate a (non unital) involutive algebra $\Ccal_c(\Tcal,A)$ which satisfies the following properties: the category of sheaf of $A$-modules over $\Tcal$ is equivalent to the category of non-degenerate $\Ccal_c(\Tcal,A)$-module, and through that equivalence, the forgetful functor from $\Ccal_c(\Tcal,a)$-module to abelian group corresponds to the functor of compactly supported section on $X$, introduced in \ref{Prop_gammac_colimit}, from sheaves of $A$-modules to abelian groups and the free $A$-module generated by $X$ corresponds to $\Ccal_c(\Tcal,A)$ seen as a module over itself.
}

\blockn{The usual case of groupoid algebra corresponds to the case where $A$ is the ring of Dedekind real or complex numbers of the topos (in the sense of \cite[D4.7]{sketches}. Assuming the axiom of dependant choice, this sheaf of ring is always admissible as soon as the topos is absolutely locally compact. Without the axiom of choice one need an additional assumption of complete regularity. But some other interesting case can be obtained for other ring, for example, Leavitt path algebras correspond to the case of a graph topos with a constant sheaf of ring (see \ref{example_Graph}).}

\blockn{Finally, as we work in a constructive context (in order to obtain a relative version of all the results) it may happens that a ringed topos satisfying all the assumptions does not have a single $A$-separating bound but a family of $A$-separating objects whose co-product is a bound. Classically a co-product of $A$-separating objects is $A$-separating so one can always take the co-product of the family as a $A$-separating bound, but constructively this is only the case when the indexing set of the co-product is decidable. To accommodate with this situation the construction of $\Ccal_c(\Tcal,A)$ is performed in the more general case of a family of $A$-separating objects instead of a single $A$-separating bound and produce a $\mathbb{Z}$-enriched pseudo-category instead of an algebra. }

\blockn{In order to obtain these results we develop under the assumption above a notion of ``compactly supported section'' of a sheaf of $A$-module $M$ over an object $X$ of the topos. When the object is ``$A$-separating'' it corresponds to the subset of sections of $M$ over $X$ which are compactly supported on $X$, but for a general $X$ compactly supported sections on $X$ are not necessarily a special case of sections and they are defined by the fact that for a fixed $M$ the compactly supported of $M$ on objects of $\Tcal$ form a cosheaf, with summation along the fiber as the functoriality.}

\blockn{Section \ref{Sec_localcompactness} contains some preliminary on locally compact Hausdorff (regular) locales and compactly supported sections (in the usual sense) over these. Section \ref{Sec_topoProp} contains preliminary on the ``topological'' assumption (separation, quasi-decidability, absolute local compactness etc...) that we will need on toposes and some of their consequences. Section \ref{Sec_Main} is the heart of the article, it contains the more general definition of compactly supported sections and all the main result of the paper. Finally, section \ref{Sec_examples} contains various example of toposes to which this construction applies and which gives back certain classical algebras.}

\blockn{We finish this introduction be some general mathematical preliminary. This paper in written in framework of constructive mathematics, we allow unbounded quantification but we do not use the unbounded replacement axiom. So it corresponds for example to the internal logic of an elementary topos as extended in \cite{shulman2010stack}. Object of this base topos $\Scal$ are called sets, and in the rest of the article by ``topos'' we mean bounded topos over $\Scal$, so ``Grothendieck toposes'' internally in $\Scal$, in the sense of categories of $\Scal$-valued sheaves over an internal site. If $\Tcal$ is a topos a  $\Tcal$ topos is a Grothendieck topos in the internal logic of $\Tcal$, in the sense that it is described by an internal site. The $2$-category of $\Tcal$-toposes is naturally equivalent to the $2$-category of toposes bounded over $\Tcal$.}

\blockn{If $\Ccal$ is a category, $|\Ccal|$ is the class of object of $\Ccal$. By pseudo category we mean a category possibly without identity elements, all the pseudo categories we will consider are enriched in abelian group, hence they are the several objects generalization of non-unital non-commutative rings.  }

\blockn{If $\Ccal$ is a site, the category of co-sheaf of sets or of abelian groups is the opposite of the category of sheaf with value in the opposite of the category of sets or of abelian group. If two sites $\Ccal$ and $\Ccal'$ defines the same topos $\Tcal$ then the category of cosheaf over $\Ccal$ and $\Ccal'$ are equivalent, this is what we call the category of cosheaf (of sets or of abelian group) over $\Tcal$.  }

\blockn{A generating family of a topos $\Tcal$ is a family of object $X_i$ such that equivalently every object can be covered by the $X_i$ or the $X_i$ form a site of definition of the topos for the induced topology. A bound of a topos is an object whose sub-object form a generating family.}

\section{Locally compact locales and compactly supported functions}

\label{Sec_localcompactness}

\blockn{Locales are the object of study of point-free topology and are very close to topological spaces. A locale $X$ is defined formally by the data of a poset $\Ocal(X)$ thought of as the poset of open subspace of $X$. More precisely, $\Ocal(X)$ must be a frame or equivalently a complete Heyting algebra, i.e. it has all supremum (called union and denoted $\cup$) and all infimum (called intersection and denoted $\cap$) and binary infimum distribute over arbitrary supremum:

\[ a \cap \bigcup_{i \in I} b_i = \bigcup_{i \in I} \left( a \cap b_i \right) \]

In particular $\Ocal(X)$ has a minimal element denoted $\emptyset$ and a maximal element denoted $X$.

A morphism of locales (also called a continuous map) $f :X \rightarrow Y$ is given by the data of an order preserving map $f^{-1}: \Ocal(Y) \rightarrow \Ocal(X)$ preserving finite infimum and arbitrary supremum. As the notation suggest, $f^{-1}$ should be thought of as the pre-image function acting on open subspaces. Every topological space $X$ defines a locale whose corresponding frame is $\Ocal(X)$, any continuous map between topological spaces defines a maps between the corresponding locales and this identifies the full subcategory of sober\footnote{a topological space in which every irreducible closed subset has a unique generic points. Every Hausdorff or locally Hausdorff topological space is sober, the underlying space of a scheme is always sober.} topological spaces with a full subcategory of the category of locales called spatial locales. For a detailed introduction to locales theory with the precise connection to topological spaces we refer the reader to \cite{borceux3} or \cite{picado2012frames}, we will use \cite{sketches} for specific results.

Any locale defines a topos of sheaves (indeed one only needs the poset of open subspaces and the notion of open coverings to define the topos of sheaves) and this construction actually identifies the category of locales with a full subcategory of the category of toposes which is called the category of localic toposes (see \cite[A4.6]{sketches}). We will not distinguishes between locales and localic toposes in the present paper.

}

\block{Locales are extremely important in topos theory, first because of the above observation, but also because for any object $X$ of a topos $\Tcal$ the poset of sub-object of $X$ form a complete Heyting algebra hence defines a locale which we will call the underlying locale of $X$ denoted\footnote{It can also be seen as the localic reflection of the slice topos $\Tcal_{/X}$ hence the notation $\loc(\Tcal_{/X})$} by $\loc(X)$ or $\loc(\Tcal_{/X})$. Any morphism $f: X \rightarrow Y $ defines naturally a continuous map $\loc(f) : \loc(X) \rightarrow \loc(Y)$ (with $f^{-1}$ being simply the pullback of sub-object along $f$.) 
}

\blockn{An example where the underlying locale of an object $X$ of $\Tcal$ appears naturally is in the description of the object $\R_{\Tcal}$ of Dedekind real numbers of $\Tcal$. This object can be described using internal logic (see \cite[D4.7]{sketches}) but it can be alternatively describe by the following universal property:

There is an isomorphism functorial in $X$:

\[ Hom_{\Tcal}(X,\R_{\Tcal}) \simeq Hom_{\loc}(\loc(X), \R) \]

where $\R$ is the locale of real numbers\footnote{Without the law of excluded middle this can be different from the topological space of Dedekind real number, see \cite[D4.7]{sketches} especially lemma $4.7.4$ and the observations after its proof.}.
Similarly, one has:

\[ Hom_{\Tcal}(X,\C_{\Tcal}) \simeq Hom_{\loc}(\loc(X), \C) \]
with $\C_{\Tcal}$ the object of Dedekind complex number and $\C$ the locale of complex number.
}

\block{Let $U,V \in \Ocal(X)$ be two open subspaces of a locale $X$, one says that $U$ is well below $V$ and write $U \ll V$ if for every directed set of open subspace $(U_i)_{i \in I}$ such that $\bigcup_{i \in I} U_i =V$ there exists $i \in I$ such that $U \leqslant U_i$. Or equivalently, if for every family $(U_i)_{i \in I}$ of open subspace of $X$ such that $\bigcup U_i = V$ one has:

\[ U \leqslant \bigcup_{j=1}^n U_{i_j} \]

for some finite family of indices $i_1,\dots, i_n$.

A locale $X$ is said to be locally compact if

\[ \forall V \in \Ocal(X), V = \bigcup_{U \ll V} U \]

A locale $X$ is said to be compact if $X \ll X$. This corresponds exactly to the usual  finite sub-covering property.
}

\block{In a locally compact locale the relation $\ll$ interpolates, i.e. if $U \ll V$ then there exists $W$ such that $U \ll W \ll V$, indeed, $V$ is the union of the $W \ll V$ and each such $W$ is the union of the $W' \ll W$ hence $V$ is the union of the $W'$ such that there exists $W$ with $W' \ll W \ll V$ such $W'$ are stable under finite union hence if $U \ll V$ then there is one of these $W'$ such that $U \leqslant W' \ll W \ll V$ which proves our claim.}

\block{One says that $U$ is rather below $V$, and denotes it $U \triangleleft V$, if there exists a $W \in \Ocal(X)$ such that $W \cup V= X$ and $W \cap U = \emptyset$. This is equivalent to the fact that the closure\footnote{The smallest closed sub-locale of $X$ containing $U$.} of $U$ is included in $V$, indeed the closed complement of $W$ corresponds to a closed sublocale of $X$ containing $U$ and included in $V$. One says that $X$ is regular if:

\[ \forall V \in \Ocal(X), V = \bigcup_{U \triangleleft V} U \]

It is a separation property essentially corresponding to the notion of regular Hausdorff (or $T3$) space.

When a locale $X$ is both locally compact and regular it is not very hard to see that:

\[ U \ll V \Leftrightarrow \left( U \triangleleft V \text{ and } U \ll X \right) \]

in more classical therms, $U \ll V$ mean that the closure of $U$ is compact and is included in $V$, $U \triangleleft V$ means that the closure of $U$ is included in $V$ and $U \ll X$ means that the closure of $U$ in $X$ is compact. These characterizations in terms of closure only works because $X$ is regular.

Moreover, for locally compact locales, being regular is equivalent to be Hausdorff (in the sense of having a closed diagonal), this is proved in \cite[II.4.8]{moerdijk2000proper}.

}

\block{\Def{A sheaf $\Fcal$ over a locally compact regular locale $X$ will be called \emph{c-soft} if for every $U \ll V$ and for every $s \in \Fcal(V)$ there exists $\widetilde{s} \in \Fcal(X)$ such that:

\[ s|_U = \widetilde{s}|_U \]

}

i.e. if when $U \ll V$ any section on $U$ which can be extended to $V$ can be extended to $X$.

This corresponds classically to the property that any section defined on a compact can be extended into a global section (because any section defined on a compact is automatically defined on a neighbourhood of the compact).
}

\block{\label{LemmeSoftExtensionInf}One says that $U \subset X$ is a neigbourhood of infinity if there exists $V$ such that $V \ll X$ and $V \cup U =X$. They corresponds to the actual neigbourhood of infinity in the one point compactification of $X$ (a presentation of one point compactification for locales is given in \cite{henry2014nonunital}).

\Lem{Let $X$ be a locally compact regular locale, $\Fcal$ a c-soft sheaf over $X$.

 let $V$ be a neighborhood of infinity and $U$ such that $U \triangleleft V$, then if $s \in \Fcal(V)$ there exists $\widetilde{s} \in \Fcal(X)$ such that:

\[ s|_U = \widetilde{s}|_U \]

}

\Dem{Assume first that $U$ is a neighborhood of infinity as well, then let $W \ll X$ such that $U \cup W =X$.

As $U \triangleleft V$ one has $U \cap W \triangleleft V$ as as $W \ll X$ one has $U \wedge W \ll X$ hence $U \cap W \ll V$. Hence as $\Fcal$ is c-soft there exists a section $s' \in \Fcal(X)$ such that $s'|_{W \cap U} = s|_{W \cap U}$. Using the sheaf condition one can then define $\widetilde{s} \in \Fcal(X)=\Fcal(U \cup W)$ to be $s'$ one $W$ and $s$ on $U$ because the two definition agrees on $W \wedge U$ and this proves the lemma in the special case.

Now in the general case, as $V$ is a neighborhood of infinity, then there exist $W$ such that $W \ll X$ and $W \cup V =X$, by interpolation one can find $W \ll W' \ll X$ and as $W \ll W'$ there is a $U'$ such that $W \wedge U' = \emptyset$ and $U' \cup W' = X$. Hence $U'$ is a neighborhood of infinity because $U' \cup W' =X$ and $U' \triangleleft V$ because $U' \wedge W = \emptyset$ and $V \cup W = X$, hence one can apply the previous special case to $U' \cup U \triangleleft V$ to obtains a section $\widetilde{s}$ which agree with $s$ on $U' \cup U$ and hence in particular on $U$.
}
}

\block{We will now discuss ``compactly supported sections''. This makes sense as soon as we are considering a sheaf $\Fcal$ which has a specific ``zero'' section marked, typically when $\Fcal$ is a sheaf of groups (that is the only case that we will consider in the present paper).

Hence let $X$ be a locally compact regular locale and $\Fcal$ be a sheaf on $X$ with a special section $0 \in \Fcal(X)$. We will say that:

\begin{itemize}

\item A section $s \in \Fcal(X)$ has support in $V \in \Ocal(X)$ if there exists $W$ such that $V \cup W =X$ and $s|_W = 0$.

\item That a section $s \in \Fcal(X)$ has compact support if it has support in $V$ for $V \ll X$.

\item That the $(U_i)_{i \in I}$ form a covering of the support of $s$ if $s$ has support in $\bigcup U_i$.

\end{itemize}

In fact one can define the support of a section $s$ has the closed complement of the open sublocale defined by ``$s=0$'' which gives the same notions as above.

One can immediately note that saying that $s$ has compact support is exactly the same as saying that $s$ is zero on some  of infinity. Also:

\Lem{ If $s$ is a compactly supported function with support in $V \in \Ocal(X)$ then there exists $V' \ll V$ such that $s$ has support in $V'$.}

\Dem{Indeed, let $W$ be a neighborhood of infinity (in the sense of \ref{LemmeSoftExtensionInf}) such that $s|_W=0$ and $W \cup V=X$ then as $V = \bigcup_{V' \ll V}  V' $ one has $W \cup \bigcup_{V' \ll V} V' =X$ if $U$ is any open subspace such that $U \cup W = X$ and $U \ll X$ there exists a $V' \ll V$ such that $U \subset V' \cup W$, in particular $V' \cup W = X$ hence $s$ has support in $V'$ and and $V' \ll V$.}

}

\block{\label{Lem_cptSup_extension}\Lem{Let $X$ be a regular locally compact locale, $\Fcal$ a c-soft sheaf of group on $X$. Let $U,V$ be two open subspaces of $X$ such that $U \ll V$ and let $s \in \Fcal(V)$. Then:

There exists a section $s' \in \Fcal(X)$ such that $s'|_U = s|_U$ and $s'$ has support in $V$. }

Using the exact same strategy as in the lemma \ref{LemmeSoftExtensionInf}, this can also be extended to the case where $U \triangleleft V$ and $V$ is a neighborhood of infinity, but as we will be mostly interested in compactly supported section in the rest of the paper we will not need this extension.

\Dem{Assume first that $U \ll V$, then let $Y$ and $Z$ such that $U \ll Y \ll Z \ll V$, let also $D$ and $D'$ be such that $Y \wedge D = Z \wedge D'=\emptyset$ and $Z \cup D = V \cup D'=X$ which exists because $Y \triangleleft Z$ and $Z \triangleleft V$. In particular as $D' \wedge Z= \emptyset$ and $Z \cup D = X$ one has $D' \triangleleft D$.

One can then take a $s_1 \in \Fcal( Y \cup D ) = \Fcal(Y) \times \Fcal(D)$ being equal to the restriction of $s$ in $\Fcal(Y)$ and $0$ one $D$. Moreover $U \cup  D' \triangleleft Y \cup D$ and $D$ is a neighborhood of infinity because $D \cup Z =X$, hence by lemma \ref{LemmeSoftExtensionInf}, there exists a section $s' \in \Fcal(X)$ which is equal to $s_1$ on $U \cup D'$, hence is equal to $s$ on $U$ and $0$ on $D'$. As $Z \cup D' =X$ hence $V \cup D' =X$ this concludes the proof in the first special case.
}
}

\block{\label{Lem_partition}\Lem{Let $X$ be a locally compact regular locale and $\Fcal$ a c-soft sheaf of abelian group over $X$. Let $s$ be a compactly supported section of $\Fcal$, and let $(U_i)$ be a covering of the support of $s$.
Then there exists $U_1,\dots, U_n$ a finite sub-cover of $U_i$ of the support of $s$ and a decomposition:

\[ s = \sum_{i=1}^n s_i \]

such that for all $i$, $s_i$ has support in $U_i$.

}

\Dem{The fact that $s$ has compact support immediately gives the finite sub-cover. We then proceed by induction on the cardinality of the covering. If $n=0$ or $n=1$ the results is obvious. Otherwise, let $W$ containing the support of $s$ and $W \ll \bigcup U_i$, then one can find $V_n$ such that $V_n \ll U_n$ and $W \leqslant V_n \cup \bigcup_{i=1}^{n-1} U_i$. Using the above lemma one can find a global section $s_n$ which agree with $s$ on $V_n$ and has support in $U_n$. Then $s-s_n$ has support included in $U_1 \cup \dots \cup U_{n-1}$ hence by induction admit a decomposition into $s-s_n=s_1 + \dots + s_{n-1}$ with each $s_i$ having support in $U_i$. Writing $s=s_1+ \dots +s_n$ concludes the proof.
}
}

\section{Absolutely separated and absolutely locally compact toposes}
\label{Sec_topoProp}

\block{An object $X$ in a topos $\Tcal$ is said to be \emph{decidable} if internally in $\Tcal$ one has $\forall x,y \in X$ $x=y$ or $x=\neq y$, externally it means that $X \times X$ can be decomposed into a disjoint sum of its diagonal and another subobject called the co-diagonal of $X$.

If $f:X \rightarrow Y$ is a map in a topos, $X$ is said to be \emph{relatively decidable} over $Y$ if $X$ is decidable as an object of $\Tcal_{/Y}$, or equivalently if $Y \times_X Y$ can be decomposed in a disjoint sum of its diagonal sub-objects and another sub-objects.
}

\block{\Def{ \begin{itemize}

\item An object $X$ in a topos is said to be \emph{quasi-decidable} if it admits a covering by a family of decidable objects.

\item A topos $\Tcal$ is said to be quasi-decidable if all its objects are quasi-decidable, i.e. if the decidable object form a generating family.

\item A geometric morphism $f: \Ecal \rightarrow \Tcal$ is said to be quasi-decidable if internally in $\Tcal$, the $\Tcal$-topos corresponding to $\Ecal$ is quasi-decidable.

\end{itemize}
}
}

\blockn{The standard terminology is to call ``locally decidable'' a topos in which every object is a quotient of a decidable object. In general, if one wants to pass from a covering by a family of decidable objects to a covering by a single decidable object it suffices to take the co-product of these decidable objects, but the co-product will be again decidable only if the indexing set is itself decidable, hence this notion is equivalent to our definition only if it is true in the ground topos that every set is a quotient of a decidable set. This induces a major difference when one apply this notion relatively: with our definition isomorphisms are quasi-decidable (any set can be covered by singletons which are decidables) in fact any localic geometric morphism is quasi-decidable, while with the more classical definition isomorphism are ``locally decidable'' exactly if it is true internally in the target that every set is a quotient a decidable set, and more generally localic morphism may fails to be ``locally decidable''. This definition of quasi-decidability is hence the correct generalization of locally decidable to a more general base topos.

\bigskip

This small distinction is not the main reason for changing the name from locally decidable to quasi-decidable. It appears that neither of these notions is actually local: the ``absolutely locally compact'' toposes we are considering in the present paper are in general only locally quasi-decidable and ``locally locally decidable'' seemed like a terminology to avoid at all cost.

}

\block{\Prop{Let $f : \Ecal \rightarrow \Tcal$ be a geometric morphism, then $f$ is quasi-decidable if and only if every object of $\Ecal$ admit a covering by an object $X \in \Ecal$ which is relatively decidable over an object of the form $f^*(I)$ for $I \in |\Ecal|$. }

\Dem{A naive external translation of the definition of quasi-decidable topos using Kripke-Joyal semantics and its extension with unbounded quantification presented in \cite{shulman2010stack} would give the following:

$f$ is quasi-decidable if and only if for all $X \in |\Tcal|$, all $v : Y \rightarrow f^* X$ in $\Ecal$, there exists an object $ S \in |\Tcal|$ such that $S \rightarrow 1_{\Tcal}$ is an epimorphism, an object $I \rightarrow S$ in $\Tcal_{/S}$ and an object $D \rightarrow f^* I$ in $\Ecal$ such that $D$ is relatively decidable over $f^* I$ and there is an epimorphism from $D$ to $Y \times f^* S$ compatible to the maps to $f^* S$.

Indeed, $X$ corresponds to the universal quantification ``for all object $Y$ of $\Ecal$'', $S$ to the existential quantification `there exists a family...'', $I$ is the indexing set of the family, $D \rightarrow f^* I$ the $I$-indexed family of decidable objects of $\Ecal$ and the epimorphism $D \rightarrow X \times f^* S$ is the covering of $X$, defined over $S$. 

It is then easy to eliminate all the redundant objects in this formulation to obtains the one in the proposition: Assuming the condition above, then taking $X=1_{\Tcal}$ and $Y$ arbitrary, one obtains $D$ relatively decidable over $f^* I$ and an epimorphism $D \twoheadrightarrow Y \times f^* S \twoheadrightarrow Y$. Conversely, assuming the condition in the proposition, let $X \in |\Tcal|$ and $v:Y \rightarrow f^*X$ in $\Ecal$, then $I \in |\Tcal|$ and $D \in |\Ecal|$ relatively decidable over $f^*I$ with an epimorphism $D \twoheadrightarrow Y$, taking $S=1_{\Tcal}$ gives all the objects of the condition above.
}

}

\block{\Cor{\begin{itemize}

\item Equivalences are quasi-decidable.

\item Localic morphisms are quasi-decidable.

\item Quasi-decidable morphisms are stable under compositions.

\end{itemize}
}

One could obviously gives internal proof of this as well, but especially for the last one the proof using the external characterization are simpler.

\Dem{\begin{itemize}
\item If $f$ is an equivalence, then every object is isomorphic to an object of the form $f^* X$.

\item If $f$ is localic then any object is a subquotient of an object of the form $f^*X$, but monomorphism $A \hookrightarrow B$ are automatically relatively decidable (their diagonal is an isomorphism and one can take the co-diagonal to be the empty) hence a subquotient of $f^* X$ is a quotient of an object relatively decidable over $f^*X$.

\item Let $f:\Fcal \rightarrow \Ecal$ and $g: \Ecal \rightarrow \Tcal$ be two quasi-decidable geometric morphisms.  Let $U$ be an object of $\Fcal$, then $U$ is a quotient of an object $U'$ relatively decidable over an object of the form $f^* V$, and $V$ is in turn a quotient of an object $V'$ relatively decidable over an object $g^* W$.

Applying $f^*$ to the objects of $\Ecal$ above one obtains that $f^* V$ is a quotient of $f^* V'$ which is relatively decidable over $(g \circ f)^* W$. Pulling back $U'$ from $f^* V$ to $f^* V'$ one obtains an object $U''$ which still cover $U$ but is now relatively decidable over $f^*V'$ which is relatively decidable over $(g\circ f)^* W$ hence $U''$ is relatively decidable over $(g \circ f)^*$ hence this concludes the proof.

\end{itemize}

}}

\block{\Prop{A pullback of a quasi-decidable geometric morphism is again quasi-decidable.}

\Dem{We work internally in the target topos. We need to show that if $\Ecal$ is quasi-decidable then $\Ecal \times \Fcal \rightarrow \Fcal$ is quasi-decidable for any topos $\Fcal$. If $\Ecal$ is quasi-decidable then it has a generating family $(D_i)_{i \in I}$ of decidable objects (take for example all the subobjects of a family of decidable objects covering a bound of $\Ecal$). The pullback of the $D_i$ to $\Ecal \times \Fcal$ are again decidable objects and they form a family of generators of $\Ecal \times \Fcal$ as a $\Fcal$ topos, hence internally in $\Fcal$, $\Ecal \times \Fcal$ has a family of decidable generators which concludes the proof.
}
}

\blockn{We also have a factorization system, that will play no role in the present paper and which we just mention for completeness:}

\block{\Prop{Let $f :\Ecal \rightarrow \Tcal$ be a geometric morphism, the full sub-category $\Ecal_{f-qd}$ of objects of $\Ecal$ which are quotients of objects relatively decidable over an object of the form $f^*X$ for $X$ an object of $\Tcal$ is stable under subobject, all colimits and all finite limit. It is in particular a topos and the inclusion $\Ecal_{f-qd} \rightarrow \Ecal$ is the $h^*$ part of a hyperconnected geometric morphism $h$ and $f$ can be factored canonically as $\Ecal \rightarrow \Ecal_{f-qd} \rightarrow \Tcal$ with the second map being quasi-decidable.}

The geometric morphisms obtained as $\Ecal \rightarrow \Ecal_{f-qd}$ are precisely those which are hyperconnected (hence $f^*$ is fully faithful and its essential image is stable under sub-objects) and such that if $X$ in $\Ecal$ is relatively decidable over $f^*Y$ then $X$ is itself $f^* Y'$. Those geometric morphism will be called anti-decidable, and this is a factorisation system.

The proof is just a routine check, and we will not use this result anyway.

}

\block{Let us now recall the definition of proper and separated morphism from \cite{moerdijk2000proper}.

\Def{\begin{itemize}
\item A topos is said to be compact if its localic reflection is compact.
\item A geometric morphism $f :\Ecal \rightarrow \Tcal$ is said to be proper if internally in $\Tcal$ the $\Tcal$-topos corresponding to $\Ecal$ is compact.
\item A Geometric morphism $f : \Ecal \rightarrow \Tcal$ is said to be separated if the diagonal morphism $\Delta_f : \Ecal \rightarrow \Ecal \times_{\Tcal} \Ecal$ is proper.
\item A Topos $\Tcal$ is said to be separated if its morphism to the terminal topos is separated, or equivalently, if the geometric morphism $\Tcal \rightarrow \Tcal \times \Tcal$ is proper.
\end{itemize}
}

Good properties of these classes of geometric morphisms (stability under pullback and composition, characterization in terms of the hyperconnected/localic factorization and so on...) are proved in \cite{moerdijk2000proper}. )

}

\block{\label{stability_prop_C'}Let $C$ be a class of maps in some category with finite limits, one says that $C$ is stable under composition if when $f \in C$ and $g \in C$ and $f \circ g $ exists then $f \circ g \in C$, and one says that $C$ is stable under pullback if when $f :X \rightarrow Y$ is in $C$ and $g: Z \rightarrow Y$ is any map then the natural projection map $\pi_2: X \times_Y Z \rightarrow Z$ is in $C$.

\Prop{Let $C$ be a class of maps which is stable under composition and pullback in a category with finite limits. Let $C'$ be the class of map $f:X \rightarrow Y $ such that the diagonal map $\Delta_f :X \rightarrow X \times_Y X$ is in $C$ then:

\begin{enumerate}
\item $C'$ is stable under composition.
\item $C'$ is stable under pullback.
\item If $f \circ g$ is in $C$ and $f$ is in $C'$ then $g$ is in $C$.
\end{enumerate}
}

\Dem{Those are all very simple ``diagram chasing'' proof. They can be found for example in \cite{moerdijk2000proper} $II.2.1$ and $II.2.2$ written in the special case where $C$ is the class of proper geometric morphisms of toposes and hence $C'$ the class of separated geometric morphism.}
}

\block{\Def{We will say that a geometric morphism $f : \Ecal \rightarrow \Tcal$ is $\Delta$-separated if the diagonal map $\Delta_f : \Ecal \rightarrow \Ecal \times_{\Tcal} \Ecal$ is a separated geometric morphism.}}

\block{\label{Prop_Deltasep}\Prop{\begin{itemize}
\item $Delta$-separated morphisms are stable under compositions and pull-back.
\item If $f \circ g $ is separated and $f$ is $\Delta$-separated then $g$ is separated.
\item If $ f \circ g$ is $\Delta$-separated then $g$ is $\Delta$-separated.
\end{itemize}}

\Dem{This follows from proposition \ref{stability_prop_C'} and the well known fact that proper geometric morphism are stable under composition and pullback (see for example \cite{moerdijk2000proper} section $I$). For the last property one also need to use proposition B3.3.8 of \cite{sketches} which say that if $f$ is any geometric morphism then $\Delta_f$ is localic, that if $f$ is localic $\Delta_f$ is an inclusion and if $f$ is an inclusion then $\Delta_f$ is an isomorphism, which implies that in particular, as any isomorphism is proper, that any geometric morphism has its diagonal map $\Delta$-separated.}
}

\block{\label{P_quasiDecImpDeltaSep}\Prop{A quasi-decidable geometric morphism is $\Delta$-separated.}

\Dem{We will work internally in the target of the morphism and show that a quasi-decidable topos is $\Delta$-separated.

Let $\Tcal$ be a quasi-decidable topos, it has a generating family $(D_i)_{i \in I}$ of decidable objects. There exists a locale $\Lcal$ whose map to the terminal locale is an open surjection and such that, internally in $\Lcal$, $I$ is covered by a sub-object of $\N$ (for example, take $\Lcal$ to be the classifying locale of the theory of partial surjections from $\mathbb{N}$ to $I$, see \cite{joyal1984extension} V.3 just after proposition $2.$ for the details ).

Hence, internally in $\Lcal$ one has a (partially) $\N$ indexed generating family of decidable object, hence internally in $\Lcal$, $\Tcal$ admit a decidable bound, hence it is localic over the Schanuel topos\footnote{The topos of set endowed with a continuous action of the localic group of permutation of $\mathbb{N}$, see \cite[C5.4.4]{sketches} } . Now the Schanuel topos is $\Delta$-separated because the localic group of automorphism of $\N$ is separated and any localic geometric morphism is $\Delta$-separated (the diagonal of a localic geometric morphism is an inclusion and inclusions are always separated because their diagonal is an isomorphism). Hence internally in $\Lcal$, the (pullback of the) topos $\Tcal$ is $\Delta$-separated. But $\Delta$-separated maps descend along open surjection because separated maps do (see \cite[C5.1.7]{sketches}), hence $\Tcal$ is $\Delta$-separated.
}

}

\block{\Def{We will say that a morphism is absolutely separated if it is separated and quasi-decidable.}

We have a work in progress whose goal is to show that the following conditions for a geometric morphism are equivalents:

\begin{itemize}

\item $f$ is absolutely separated.

\item $f$ is separated and $\Delta$-separated.

\item $f$ is strongly separated, i.e. the diagonal of $f$ is tidy (cf. \cite{moerdijk2000proper}).

\end{itemize}

The equivalence of the last two conditions is already in \cite{moerdijk2000proper}, and we proved above that they are implied by the first condition. We hope to be able to publish a proof of the converse implication soon.

}

\blockn{We recall the main result of \cite{henry2015finitness}:}

\block{\label{HCfinitnessTh}\Th{A Topos which is hyperconnected and absolutely separated has a generating familly of objects which are internally finite and decidable. }

Moreover this result is valid within the internal logic of a topos.

}

\block{\Def{An object $X$ of a topos $\Tcal$ is said to be quasi-finite of cardinal $n$ object if internally in $\Tcal$ one has $\exists x \in X \Rightarrow  X$ is finite decidable of cardinal $n$.
One says that it is quasi-finite if it is quasi-finite of cardinal $n$ for some $n$.
}

I.e. quasi-finite objects are not necessarily finite, but there is a sub-terminal object called there support such that they are finite of cardinal $n$ over their support. The support is just the interpretation of the proposition $\exists x \in X$.
}

\block{\label{PropQfinExist}\Prop{An absolutely separated topos has a basis of quasi finite objects.}

\Dem{If one has an absolutely separated topos $\Tcal$, then the hyperconnected geometric morphism to its localic reflection $\Lcal$ is also quasi-decidable by \cite[5.2]{henry2015finitness} and separated by \cite[II.2.5]{moerdijk2000proper}, hence one can apply theorem \ref{HCfinitnessTh} above to it and obtains internally in $\Lcal$ a basis of finite object. Local sections of the stack over $\Lcal$ of finite objects of $\Tcal$ are objects of $\Tcal$ that are finite over their domain of definition, one can then restrict to the open subspace of $\Lcal$ where their cardinal is equal to $n$ for any fixed $n$ and one get quasi-finite objects of $\Tcal$, and those clearly form a generating family. }

}

\block{\label{def_absloccpt}\Def{\begin{itemize}
\item A topos $\Tcal$ is said to be absolutely compact if it is compact and absolutely separated.

\item A topos $\Tcal$ is said to be absolutely locally compact if it has a generating family of objects $X$, such that the topos $\Tcal_{/X}$ is absolutely separated with a locally compact localic reflection.

\item An object $X$ of an absolutely locally compact topos $\Tcal$ is said to be separating if $\Tcal_{/X}$ is absolutely separated and has a locally compact localic reflection. 

\end{itemize}
}}

\block{Under the conjecture mentioned below that absolutely separated is equivalent to strongly separated, absolutely compact would mean that the topos, its diagonal, and the diagonal of the diagonal are all proper. We also conjecture that a topos is absolutely locally compact if and only it is exponentiable, its diagonal is exponentiable and the diagonal of its diagonal is also exponentiable.}

\block{\label{Lem_ll_ext_int}\Lem{Let $\Tcal$ be a separated topos, let $X$ be an object of $\Tcal$ such that $\loc(\Tcal_{/X})$ is locally compact and Hausdorff. Let also $U$ be a subobject of $X$ such that $U \ll X$ (in $\loc(\Tcal_{/X})$). Then internally in $\Tcal$ there exists a finite object $F$ such that $U \subset F \subset X$. }

We insist on the fact that we do not claim that such an object exists externally, this is not true in most situation.

\Dem{$U$ can be seen as an open subspace of $\loc(\Tcal_{/X})$. Let $\overline{U}$ the closure of $U$ in this locale, as the locale is locally compact and Hausdorff and $U \ll X$ the closure is compact. We denote by $\Tcal_{/\overline{U}}$ the pullback of $\Tcal_{/X} \rightarrow \loc(\Tcal_{/X})$ along $\overline{U} \rightarrow \loc(\Tcal_{/X})$. By \cite[C2.4.12]{sketches} localic morphisms and hyperconnected morphisms are stable under pullback, hence the morphism $\Tcal_{/\overline{U}} \rightarrow \Tcal_{/X}$ is localic and the morphism $\Tcal_{/\overline{U}} \rightarrow \overline{U}$ is hyperconected, in particular, $\overline{U}$ is the localic reflection of $\Tcal_{/\overline{U}}$ hence (as $\overline{U}$ is compact) $\Tcal_{/\overline{U}}$ is proper as a topos and hence, as $\Tcal$ is separated the geometric morphism $\Tcal_{/\overline{U}} \rightarrow \Tcal$ is proper, as it is localic it means that it corresponds to a compact locale in the internal logic of $\Tcal$. Internally in $\Tcal$ one has hence an inclusion $U \subset \overline{U} \subset X$ with $U$ and $X$ two sets and $\overline{U}$ a compact sub-locale of $X$ (seen as a discrete locale). But (still internally in $\Tcal$) $X$ can be covered by the $\{ x \}$ for $x \in X$, hence there exists a finite collection of them which cover $\overline{U}$, and in particular $U$, the union of such a finite collection produces (internally) the desired finite set. }}

\block{\label{P_DeltaSepvsAbslocSep} The separating objects of an absolutely locally compact topos are not necessarily decidable. In fact:

\Prop{Let $\Tcal$ be an absolutely locally compact topos, then the following conditions are equivalent:

\begin{enumerate}

\item $\Tcal$ is quasi-decidable.

\item $\Tcal$ is $\Delta$-separated.

\item Every separating object of $\Tcal$ is decidable.

\item $\Tcal$ admit a generating family of decidable separating objects.

\end{enumerate}

}

\Dem{The implication $1. \Rightarrow 2.$ has been proved for a general topos in \ref{P_quasiDecImpDeltaSep}.

For the implication $2. \Rightarrow 3.$, if $X$ is a separating object of $\Tcal$ then $\Tcal_{/X}$ is separated, hence if $\Tcal$ is $\Delta$-separated the morphism from $\Tcal_{/X}$ to $\Tcal$ is separated by the third point of \ref{stability_prop_C'} and by \cite[II.1.3]{moerdijk2000proper} the map $\Tcal_{/X} \rightarrow \Tcal$ is separated precisely when $X$ is decidable.

The implications $3. \Rightarrow 4.$ and $4. \Rightarrow 1.$ are trivial. }

}

\block{We conclude this section with an important example of the above proposition:

\Prop{Let $\Tcal$ be an etendu, i.e. the topos of equivariant sheaves on an étale localic groupoid $G_1 \rightrightarrows G_0$.

Assume additionally that $G_0$ is separated (Hausdorff), and hence that $G_1$ is locally separated then the following condition are equivalent:

\begin{enumerate}

\item $\Tcal$ satisfies the equivalent condition of proposition \ref{P_DeltaSepvsAbslocSep} above.
\item $G_1$ is separated (Hausdorff).

\end{enumerate}

}

In the next section, we will see that the construction of the convolution algebra of a topos will be more subtle in the case where the separating objects are not decidable. by this proposition this corresponds (in the case of an etendu) to the situation of a non-Hausdorff étale groupoid. It will be rather evident for a reader familiar to the construction of the convolution algebra of an étale groupoid that this subtleties are exactly similar to the fact when we construct the $C^*$-algebra of an etale groupoid we need to replace continuous function on $G_1$ by linear combination of functions which are compactly supported on Hausdorff open subsets of $G_1$ and extended by $0$ on the rest of $G_1$.

\Dem{We will check that condition $2.$ is equivalent to the fact that $\Tcal$ is $\Delta$-separated which is one of the conditions of proposition \ref{P_DeltaSepvsAbslocSep}.

The fact that $\Tcal$ is the topos of the groupoid $G_1 \rightrightarrows G_0$ can be translated in the fact that one has a pullback square:

\[\begin{tikzcd}[ampersand replacement = \&]
G_1 \arrow{r} \arrow{d} \& G_0 \arrow{d} \\
G_0 \arrow{r} \& \Tcal \\
\end{tikzcd} \]

This square gives rise to another pullback square:

\[\begin{tikzcd}[ampersand replacement = \&]
G_1 \arrow{r} \arrow{d} \& \Tcal  \arrow{d} \\
G_0 \times G_0 \arrow{r} \& \Tcal \times \Tcal \\
\end{tikzcd} \]

So as $G_0 \times G_0 \rightarrow \Tcal \times \Tcal$ is an etale covering, $\Tcal$ is $\Delta$-separated if and only if $G_1 \rightarrow G_0 \times G_0$ is separated, but as $G_0 \times G_0$ is separated and $\Delta$-separated (every locale is $\Delta$-separated) this is equivalent to the fact that $G_1$ is separated by \ref{stability_prop_C'}.

}
}

\section{Admissible sheaf of rings and the algebra of compactly supported functions }

\label{Sec_Main}

\block{\label{DefAdmissibleSR}\Def{Let $\Tcal$ be an absolutely locally compact topos and let $A$ be a (possibly non-commutative) unital ring object over $\Tcal$. One says that $A$ is admissible if there is a generating family $(X_i)_{i \in I}$ of separating objects of $\Tcal$ such that for each $i \in I$:

\begin{itemize}

\item When seen as a sheaf on the locale $Loc(\Tcal_{/X_i})$ by restriction to the sub-object of $X_i$, $A$ is c-soft.

\item $\Tcal_{/X_i}$ admit a generating family of quasi-finite objects of cardinals invertible in $A$.

\end{itemize}
}

Separating objects $X_i$ such that these two conditions holds will be called $A$-separating objects of $\Tcal$. Of course, if every integer is invertible in $A$ the second condition always holds.

The key example we have in mind is when $A$ is the sheaf of Dedekind real numbers or complex numbers, or eventually some algebra over this, for example an internal $C^*$-algebra or internal Banach algebra:
}

\block{\Prop{Assuming the axiom of dependant choice, for every absolutely locally compact topos if the ring $A$ contains the ring of Dedekind real numbers then $A$ is admissible and every separating object is $A$-separating. }

If we do not assume the axiom of dependant choice one needs to assume additionally that there is a generating family of separating objects $X_i$ such that the localic reflection of $\Tcal_{/X_i}$ is completely regular in addition of being locally compact and Hausdorff, which become automatic by Urysohn's lemma if one assumes the axiom of dependant choice.

\Dem{In this case, every integer is invertible in $A$ hence the second condition of definition \ref{DefAdmissibleSR} is automatic because of \ref{PropQfinExist}. The first conditions follow easily from complete regularity (which comes from Urysohn lemma, see \cite[XIV.5 and 6]{picado2012frames}) that allows to construct for any $U \ll V$ two subobjects of $X$ a positive real valued functions $f$ which is $1$ on $U$ and zero outside of $V$. then any section defined on $V$ can be extended to $X$ by multiplying by $f$ and extending by $0$ outside of $V$. }
}

\block{\label{Lem_partition2} If $X$ is a separating object of a topos $\Tcal$, or more generally any object such that $\Tcal_{/X}$ is locally compact, an arrow $s:X \rightarrow M$ for $M$ some abelian group object in $\Tcal$ can be said to be compactly supported on $X$ if it is compactly supported when seen as a section in $Loc(\Tcal_{/X})$ of the restriction of $M$ the $Loc(\Tcal_{/X})$, or to put it more explicitly if $X=U \cup W$ with $U \ll X$ and $s|_W=0$.

\Lem{Let $\Tcal$ be an absolutely locally compact topos, $A$ an admissible ring object over $\Tcal$ and $M$ an $A$-module in $\Tcal$ and $X$ a $A$-separating object.

Let $m$ be a compactly supported section of $M$ over $X$ and let $(U_i \hookrightarrow X)_{i=1 \dots n}$ a covering of the support of $m$ by sub-object of $X$. Then there exists $(\lambda_i)_{i=1 \dots n}$ compactly supported sections of $A$ on $X$ such that:

\begin{itemize}

\item for all $i$, $\lambda_i$ has support in $U_i$,

\item $\displaystyle \sum_{i=1}^n \lambda_i m =m. $

\end{itemize}

}
\Dem{Let $V$ be a sub-object of $X$ such that $V \ll \bigcup U_i$ and $V$ contains the support of $m$. Pick any $W$ such that $V \ll W \ll X$, as $A$ is c-soft over $X$, one has by lemma \ref{Lem_cptSup_extension} a section $\lambda \in A(X)$ such that $\lambda|_V=1$ and $\lambda$ has support in $W$ (in particular $\lambda$ is compactly supported),
and by lemma \ref{Lem_partition} one can find $\lambda_1,\dots,\lambda_n$ such that $\lambda_i$ has support in $W \wedge U_i$ and $\sum \lambda_i = \lambda$.

As $\lambda$ is equal to $1$ on $V$ and $m$ has support in $V$ one has $\lambda m =m$ and hence $\sum \lambda_i m = m$ which concludes the proof. 
}
}

\block{\label{Prop_gammac_colimit}Let $\Tcal$ be an absolutely locally compact topos, and $A$ be an admissible ring object over $\Tcal$. 
Let $V$ be a sheaf of $A$-module over $\Tcal$.

For any $A$-separating object $X$ of $\Tcal$ we define $\Gamma_c(X; V)$ to be the set of sections of $V$ over $X$ which are compactly supported.

If $f : V \rightarrow W$ is a map between two sheaves of $A$-modules, then post-composing with $f$ induces a map $f :\Gamma_c(X,V) \rightarrow \Gamma_c(X,W)$, which makes $\Gamma_c(X, \_ )$ into a functor.

\Prop{Let $A$ be an admissible sheaf of ring over an absolutely locally compact topos. Then for any $A$-separating object $X$, the functor $\Gamma_c(X, \_ )$ from the category of $A$-modules in $\Tcal$ to the category of commutative groups commute to all colimits and all finite limits.}

\Dem{

\begin{itemize}

\item As $\Gamma_c(X, \_)$ is an additive functor between additive categories, it commutes to the zero object and to bi-products, hence both to finite co-products and finite products.

\item $\Gamma_c(X, \_)$ commute to kernel (and hence to all finite limite): let $f:V \rightarrow W$ be a morphism of sheaf of $A$-module and let $P \hookrightarrow V$ be the kernel of $f$. A map with compact support from $X$ to $P$ is exactly a map $a$ from $X$ to $V$ which has compact support and such that $f \circ a = 0$, hence $\Gamma_c(X,P)$ is exactly the kernel of $f: \Gamma_c(X,V) \rightarrow \Gamma_c(X,W)$.

\item $\Gamma_c(X, \_)$ commutes to all finite colimits: as we already know that it is a left exact functor between abelian categories it is enough, for example by \cite[1.11.2 and 1.11.4]{borceux2}, to show that it send epimorphisms to epimorphisms.

Let $f:  V \twoheadrightarrow W$ be an epimorphism between two sheaves of $A$-modules, and let $s: X \rightarrow W$ be a map with compact support, our goal is to lift it into a compactly supported function from $X$ to $V$.
As $f$ is an epimorphisms, there exists a covering $v_i:X_i \rightarrow X$ with for each $i$ a map $s_i : X_i \rightarrow V$ such that $f \circ s_i = s \circ v_i$. Because $X$ is $A$-separating one can assume that for each $i$, $X_i$ is quasi-finite of cardinal $k_i$, that its support is relatively compact in $X$ and that $k_i$ is invertible in $A$ and moreover there is a finite family $X_1, \dots, X_n$ of such object that cover the support of $s$.

One can then (be lemma \ref{Lem_partition2}) find functions $\lambda_1, \dots , \lambda_n$ from $X$ to $A$ such that for each $i$, $\lambda_i$ is supported in the support of $X_i$ and:

\[ \sum_{i=1}^n \lambda_i s = s \]

Finally, we define internally in $\Tcal_{/X}$:

\[ s' := \sum_{i = 1}^ n \left\lbrace
\begin{array}{ccc}
\frac{\lambda_i}{k_i} \displaystyle \sum_{v \in X_i } s_i(v)  & \mbox{if} & \exists x \in X_i\\
0 & \mbox{if} & \lambda_i =0 \\
\end{array}\right.\]

It is well defined: for each $i$, as $\lambda_i$ is supported in the support of $X_i$ hence internally in $\Tcal_{/X}$ one has $\exists x \in X_i$ or $\lambda_i=0$. In the case $\exists x \in X_i$, the object $X_i$ of $\Tcal_{/X}$ is finite and decidable hence the summation over it is well defined and if the two condition holds simultaneously they both gives the results $0$.

$s'$ is a map from $X$ to $V$, its support is included in the union of the support of the $X_i$, hence it is a compactly supported map from $X$ to $V$ and:

\[ f(s') := \sum_{i = 1}^n \left\lbrace
\begin{array}{ccc}
\sum_{v \in X_i } \frac{\lambda_i s}{k_i}  & \mbox{if} & \exists x \in X_i\\
0 & \mbox{if} & \lambda_i =0 \\
\end{array}\right.
 = \sum_{i=1}^n \lambda_i s =s
\] 

hence $f : \Gamma_c(X,V) \rightarrow \Gamma_c(X,W)$ is also an epimorphism.

\item $\Gamma_c(X, \_)$ commutes to arbitrary co-product: In the case of decidable co-product there is a very simple proof, any compactly supported map to the coproduct will factor through a finite (decidable) co-product by compactness and we already have commutation to finite co-product so the proof is finished. But unfortunately, this argument is not sufficient for a non decidable coproduct and unless the base topos is assumed to be locally decidable, decidable coproduct will not suffice. For the general case we will need the lemma \ref{Lem_gammac_indexedcoprod} below.

Let $(V_i)_{i \in I}$ be a family of $A$-module.The comparison map:

\[ \bigoplus \Gamma_c(X,V_i) \rightarrow \Gamma_c \left( X, \bigoplus_i V_i \right) \]

is always a monomorphism because we already know that $\Gamma_c(X,\_)$ preserve monomorphism hence each component is a monomorphism, so it is enough to show that it is an epimorphism, i.e. that each element of $\Gamma_c \left( X, \bigoplus V_i \right)$ can be written as a finite combination of elements in the $\Gamma_c(X,V_i)$.

The family $(V_i)_{i \in I}$ can be seen as a $A$-module $V$ over $Y=p^* I$, whose internal direct sum $W$ is just the ordinary direct sum. Hence, any compactly supported map $s$ from $X$ to $W$ can be described, as in lemma \ref{Lem_gammac_indexedcoprod} below, by an object $D$ with a map $D \rightarrow f^* I$, a map $D \rightarrow X$ and a compactly supported section $s'$ from $D$ to $V$ in $\Tcal_{/f^*I}$.

Such a map $D \rightarrow f^* I$ can be interpreted as a decomposition of $D$ into $D =\coprod D_i$. The function $s' :D \rightarrow \coprod V_i$ is hence a collection of compactly supported functions $s_i : D_i \rightarrow V_i$, which are all zero except a finite number. Using the formula in the lemma those functions $s_i$ can be then turned back into compactly supported functions from $X$ to $V_i$ whose sum is the initial function $s$ from $X$ to $W$ and this concludes the proof.

\end{itemize}

}
}

\block{\label{Lem_gammac_indexedcoprod}\Lem{Let $X$ be a $A$-separating object of $\Tcal$ and let $Y$ be any object, let $V$ be a sheaf of $A$-modules in $\Tcal_{/Y}$, let $W$ be the sheaf of $A$-modules defined internally in $\Tcal$ by:

\[ W= \bigoplus_{y \in Y} V_y \]

Let also $s:X \rightarrow W$ be a compactly supported function.

Then there exists a $A$-separating object $D$ with two maps $p_1 :D \rightarrow X$ and $p_2: D \rightarrow Y$, and a compactly supported function $s'$ from $(D,p_2)$ to $V$ in $\Tcal_{/Y}$ such that (internally) for all $x \in X$:

\[ s(x) = \sum_{p_1(d)=x} s'(d)  \]

}


\bigskip

This lemma is key for several results in the paper. It roughly says that compactly supported function to an internal co-product can be written in some sense as ``compactly indexed sum'' of section of the component.

\Dem{We will first explain why the formula for $s(x)$ is meaningful. As $D$ and $X$ are separating, any map from $D$ to $X$ is fiberwise decidable, and there is a subobject $D' \subset D$ such that $D'$ is finite over $X$ and contains the support of $s'$, hence the sum can be seen (internally) as a finite sum. Moreover $s'(d)$ is an element of $V_{p_2(d)}$, but as $W= \bigoplus_{y\in Y} V_y$, it can be seen as an element of $W$ too.

We will now prove the lemma. Let $X,Y,V,W,s$ be as in the lemma. By definition of $W$, internally in $\Tcal$, for all $x \in X$ , there exists an integer $n$, element $y_1, \dots ,y_n \in Y$ and $v_1 \in V_{y_1}, \dots, v_n \in V_{y_n}$ such that $s(x) = \sum v_{y_i}$.

This internal statement can be attested by a covering $v_i: X_i \rightarrow X$ of $X$, for each $i$, an integer $n_i$, $n_i$-applications $y^i_1, \dots y^i_{n_i} : X_i \rightarrow Y$, $n_i$-applications $v^i_j: X_i \rightarrow V$ over $Y$ such that:

\[ \sum_{j=1}^{n_i} \iota(v^i_j) = s \circ v_i \]

where $\iota$ denote the canonical arrow from $V$ to $W$.

Moreover, because $X$ is $A$-separating one can freely assume that each $X_i$ is quasi-finite of cardinal $m_i$ over $X$, and that image of $X_i$ in $X$ is relatively compact in $X$. Moreover one can extract a finite family $X_1,\dots,X_k$ which covers the support of $s$, as well as (by lemma \ref{Lem_partition2}) a family of functions $\chi_1,\dots, \chi_k : X \rightarrow A$ such that $\chi_i$ has its support contains in the support of $X_i$, and $\sum \chi_i s =s $, in fact we can (and we will) further assume that $\chi_i$ is compactly supported within the support of $X_i$.

Let then:

\[ D =\coprod_{i=1}^k \coprod_{j=1}^{n_i} X_i \]

$D$ is a decidable coproduct of separating object hence it is a separating object. Let $p_1 : D \rightarrow X$ be the natural map that send each $X_i$ to $X$ by $v_i$. Let $p_2 :D \rightarrow Y$ be the map that send $(i,j)$ component $X_i$ to $Y$ by $y^i_j$ and $s': D \rightarrow V $ be the map that send the $(i,j)$ component $X_i$ of $D$ to $V$ by $ \frac{v_i^j \chi_i}{m_i}$

As the $\chi$ are compactly supported within the image of $X_i$ in $X$ and $X_i$ is finite over its image, this function from $D$ to $V$ is indeed compactly supported and:

\begin{multline*}  \sum_{p_1(d)= x} s'(d) = \sum_{i=1}^k \sum_{j=1}^{n_i} \sum_{v_i(d)=x} \frac{\chi_i(x) v_i^j(d)}{m_i} =  \sum_{i=1}^k\sum_{v_i(d)=x} \frac{\chi_i(x)}{m_i}  \sum_{j=1}^{n_i}  v_i^j(d) \\ =  \sum_{i=1}^k\sum_{v_i(d)=x} \frac{\chi_i(x)}{m_i} s \circ v_i(d) = \sum_{i=1}^k \chi_i(x) s(x) = s(x)  \end{multline*}

} 
}

\block{\label{Gamma_functo_left}We will now describe the functoriality in $X$ of $\Gamma_c(X;V)$. Let $f : X \rightarrow Y$ be an arrow between two $A$-separating objects of an absolutely locally compact topos $\Tcal$. As $\Tcal_{/X}$ is separated and $\Tcal_{/Y}$ is $\Delta$-separated, the map $\Tcal_{/X} \rightarrow \Tcal_{/Y}$ is separated by \ref{stability_prop_C'} and hence $(X,f)$ is a decidable object of $\Tcal_{/Y}$ (see for example \cite[II.1.3]{moerdijk2000proper}). Let $v:X \rightarrow V$ be a compactly supported section on $X$ of some abelian group object $V$. There exists a sub-object $U \ll X$ and $W \subset X$ such that $U \cup W = X$ and $v|_W=0$, by lemma \ref{Lem_ll_ext_int}, internally in $\Tcal_{/Y}$ there exists a finite object $F$ such that $U \subset F \subset X$, in particular internally in $\Tcal_{/Y}$ one can define:

\[ \sum_{x \in X} v(x)  \]

as the sum for $x \in F$ because $F$ is finite and decidable and contains the support of $v$, and this does not depend on the  (internal) choice of $F$. In fact one has:

\Prop{Let $\Tcal$ be an absolutely locally compact topos, $f :X \rightarrow Y$ be an arrow between two separating objects of $\Tcal$. And $v:X \rightarrow V$ be a compactly supported arrow to an abelian group. Then internally in $\Tcal$ the following formula:

\[ w(y) := \sum_{f(x) = y} v(x) \]

defined a compactly supported function from $Y$ to $V$.
}

$w$ will be denoted $\Sigma_f v$ and this turns $\Gamma_c( \_ , \_)$ into a bi-functor.

\Dem{In the discussion above, we proved that\footnote{proving something internally in $\Tcal_{/Y}$ is exactly the same as proving that internally in $\Tcal$ the same thing holds for all $y \in Y$.} internally in $\Tcal$, for each $y \in Y$ there exists a finite set $F \subset f^{-1}(\{y\})$ such that for all $x \in f^{-1}(\{y\})$ either $v(x)=0$ or $x \in F$ which proves that the above sum is well defined and defines a function from $Y$ to $V$. We just have to prove that it is compactly supported, but the set of $f^*(V)$ such that $V \ll Y$ form a directed covering of $X$, hence as $v$ is compactly supported there exists a  $V \ll Y$ such that $U$ (the ``support'' of $v$) is included in $f^*(V)$ . Then, for all $y \in Y$ let $F$ be a finite set as above, then either $v=0$ at every element of $F$, in which case $w(y)=0$ or there exist an element of $F$ which is in $U$, in which case $y \in V$, hence for all $y \in Y$, $w(y)=0$ or $y \in V$ with $V \ll Y$ externally, hence $w$ is compactly supported which concludes the proof. }
}

\block{\label{Lem_WeakseparatingTransport}
Let $\Tcal$ be an absolutely locally compact topos in which the terminal object $1$ is separating. Let also $X$ be a separating object of $\Tcal$, in particular, $\Tcal_{/X}$ is separated and as $1$ is separating, the topos $\Tcal = \Tcal_{/1}$ is $\Delta$-separated hence $\Tcal_{/X} \rightarrow \Tcal$ is separated (by \ref{Prop_Deltasep}) which means that $X$ is decidable (by \cite[II.1.3]{moerdijk2000proper}). For any abelian group object $V$ of $\Tcal$ and every object $X$ of $ \Tcal$ one can define the abelian group object $\bigoplus_{x \in X} V$, but in the special case where $X$ is decidable it can be identified with the group of finitely supported functions from $X$ to $V$.

In particular both $\Gamma_c(X,V)$ and $\Gamma_c(1, \bigoplus_{x \in X} V)$ corresponds to subgroup of the group of all functions from $X$ to $V$.

\Lem{In the situation above, $\Gamma_c(X,V)$ and $\Gamma_c(1, \bigoplus_{x \in X} V)$ are equal as subgroup of $Hom(X,V)$, moreover this identification of $\Gamma_c(X,V)$ and $\Gamma_c(1, \bigoplus_{x \in X} V)$ is functorial in both $X$ and $V$.}

\Dem{Let $v: X \rightarrow V$ be a compactly supported function and let $U \ll X$ which contains the support of this function. By lemma \ref{Lem_ll_ext_int}, internally in $\Tcal$ there exists a finite object $F$ such that $U \subset F \subset X$, i.e. $v$ is a finitely supported function from $X$ to $V$ and hence corresponds to a map $1 \rightarrow \bigoplus_{x \in X} V$ it is compactly supported by the exact same argument as in the end of the proof in \ref{Gamma_functo_left}.

Conversely, if $v:1 \rightarrow \bigoplus_{x \in X} V$ is a compactly supported functions, then one can apply lemma \ref{Lem_gammac_indexedcoprod}, and one obtains an object $D$ with a compactly supported map $\lambda:D \rightarrow V$, a map $p_2:D \rightarrow X$ and $p_1$ the map $p_1:D \rightarrow 1$ such that:

\[ v = \sum_{d \in D} i_{p_2(d)}(\lambda(d))  \]

where for $x \in X$, $i_{x}$ denotes the corresponding map $V \rightarrow \bigoplus_{x \in X} V$. In particular seeing $v$ as a function from $X$ to $V$ one has exactly:

\[ v(x) = \sum_{p_2(d)=x} \lambda(d) \]

hence $v= \Sigma_{p_2} \lambda$ is indeed an element of $\Gamma_c(X;V)$ by \ref{Gamma_functo_left}. 
The fact that this identification is functorial is immediate. 
}}

\block{\label{Th_def_cptsection}\Th{Let $\Tcal$ be an absolutely locally compact topos, $A$ an admissible ring object of $\Tcal$. Then one has a unique bifunctor $\Gamma_c(X;V)$ for $X \in |\Tcal|$ and $V$ a sheaf of $A$-module on $\Tcal$ such that:

\begin{itemize}

\item For all $V$, $\Gamma_c(\_ ,V)$ is a cosheaf of abelian groups on $\Tcal$.

\item When $X$ is $A$-separating this coincide with the definition in \ref{Prop_gammac_colimit}.

\end{itemize}

Moreover, for all $X$, $\Gamma_c(X,\_)$ commutes to all colimits.

}

This theorem will be our definition of a ``compactly supported section of $V$ on $X$'' when $X$ is not separating: they are the elements of $\Gamma_c(X;V)$.

\Dem{For the first part of the proposition, it is enough to proves that for any $V$ a sheaf of $A$-module, $\Gamma_c(X;V)$ as defined in \ref{Prop_gammac_colimit} for $X$ a $A$-separating object defines a cosheaf of abelian group for the canonical topology of $\Tcal$. Indeed the category of $A$-separating objects endowed with the canonical topology of $\Tcal$ is a site of definition for $\Tcal$, so this will construct a cosheaf of abelian groups $\Gamma_c( \_ ;V )$ on $\Tcal$ for all $V$.

Let $X$ be a separating object, and let $U_i$ be a covering of $X$ by separating object, we need to prove the cosheaf condition, which can be formulated as follow: $X$ can be written as a certain co-limits of the $U_i$ and there fiber products, and one needs to show that $\Gamma_c(\_,V)$ commutes to this colimit. As the colimit is computed in $\Tcal_{/X}$ one can freely assume that $X=1$ by working in $\Tcal_{/X}$. But then $\Gamma_c(Y,V) = \Gamma_c(1, \bigoplus_{y \in Y} V)$ obviously commutes to all co-limits: $Y \mapsto \bigoplus_{y \in Y} V$ commutes to all co-limit because it is the left adjoint to the functor which send a sheaf of $A$-module $W$ to the sheaf (of sets) of morphism from $V$ to $W$ and $\Gamma_c(1,\_)$ commutes to all co-limits because of proposition \ref{Prop_gammac_colimit}.

For the last claim of the proposition, co-limits in the category of co-sheaves over a site are computed objectwise so the fact that $\Gamma_c(X; \_)$ commutes to all co-limit when $X$ is separating shows that the functor $V \mapsto \Gamma_c( \_;V)$ commutes to co-limits as a functor from sheaves of $A$-module to co-sheaf of abelian groups, and hence that for all $X$, $\Gamma_c(X,\_)$ commutes to all co-limits.
}
}

\block{While the definition of $\Gamma_c(X;V)$ for a general $X$ may sound very abstract, it is not hard to give explicit formulas to compute them using the cosheaf property: Let $D$ be any separating object covering $X$ and let $D'$ be a separating object covering $D \times_X D$. Then $X$ is the coequalizer in $\Tcal$ of the two maps from $D$ to $D'$, hence $\Gamma_c(X;V)$ is the co-equalier of $\Gamma_c(X;D') \rightrightarrows \Gamma_c(X;D)$. If the topos is quasi-decidable $D \times_X D$ will itself be separating and can be used instead of $D'$.

The above description works only if we are able to find a single separating object covering an object $X$, which is the case as soon as the base topos is boolean or if in the base topos every object can be covered by a decidable object. If it is not the case one has the following more general description: pick $D_i$ a covering family of $X$ by separating object, and for all $i,j$ pick\footnote{This can always be done without invoking the axiom of choice by using the collection axiom introduced in \cite{shulman2010stack}.} $D'_{i,j,k}$ a covering family of $D_i \times_X D_j$, then $\Gamma_c(X;V)$ can be computed as the coequalizer:

\[ \bigoplus_{i,j,k} \Gamma_c(D'_{i,j,k};V) \rightrightarrows \bigoplus_i \Gamma_c(D_i;V) \twoheadrightarrow \Gamma_c(X;V)  \]

}

\block{\label{Prop_Gammac_exchange}\Prop{If $\Tcal$ is a absolutely locally compact topos with $A$ an admissible sheaf of ring, then one has an isomorphism functorial in $V,X$ and $Y$:

\[ \Gamma_c( Y \times X;V) \simeq \Gamma_c \left( Y,\bigoplus_{x \in X} V \right) \]

}
\Dem{When $Y=1$ and both $1$ and $X$ are $A$-separating, this is lemma \ref{Lem_ll_ext_int}. Assuming $1$ is $A$-separating, then the two side defines co-sheaves of abelian groups in $X$ and are functorially isomorphic when $X$ is separating and one has such an isomorphism for all $X$.

In particular, for a general absolutely locally compact topos $\Tcal$, for any $A$-separating object $Y$ and any object $X$ one has an isomorphism:

\[ \Gamma_c(Y \times X, V) \simeq \Gamma_c(Y, \bigoplus_{x \in X} V) \]

by applying the above result in $\Tcal_{/Y}$ to the object $X\times Y$. But here again, the two sides are cosheaves in $Y$ hence the isomorphism extend to all $Y$.
}  
}

\block{\label{Prop_MatrixRep}Let $X$ be an object of $\Tcal$, we denote by $X_A$ the free $A$-module generated by $X$, i.e., internally in $\Tcal$:

\[ \bigoplus_{x \in X} A =X_A \]

A morphism from $X_A$ to any other $A$-module $M$ is the same as a map from $X$ to $M$. If $X$ is separating, it hence makes sense to ask whether such a map $X_A \rightarrow M$ is compactly supported (depending on if the corresponding map $X \rightarrow M$ is compactly supported or not. The result above shows that:

\Th{Let $X,Y$ be two $A$-separating objects of an absolutely locally compact topos with an admissible ring object $A$. Then a compactly supported map $X_A \rightarrow Y_A$ is the same as a compactly supported section $\Gamma_c(X \times Y; A)$.

\bigskip

The correspondence between the two is as follow:

\bigskip

If $\gamma \in \Gamma_c(X \times Y;A)$ and if one has $p=(p_1,p_2):D \rightarrow X \times Y$ with $D$ a separating object such that $\gamma = \sigma_p \gamma_0$ then the map $F: X \rightarrow Y_A$ corresponding to $\gamma$ is described internally as:

\[ F(x) = \sum_{p_1(d)=x} i_{p_2(d)}(\gamma_0(d)) \]

Where the $i_y$ for $y \in Y$ are the structural maps from $A$ to $Y_A$ 

}

This theorem is one of the key result of this paper. It should be understood as a description of compactly supported map from $X_A$ to $Y_A$ by matrix elements $X \times Y \rightarrow A$, but with the subtleties that matrix elements are not a function from $X \times Y$ to $A$ as one should expect, but a compactly supported section, and that in the case where $X \times Y$ is not decidable, those compactly supported sections are not ``sections which are compactly supported''.

\Dem{The isomorphism of \ref{Prop_Gammac_exchange} gives us directly that:

\[ \Gamma_c( X;A Y) \simeq \Gamma_c(X \times Y; A) \]

We just need to show that it is indeed as described in the theorem, which amount to understand the composite:

\[ \Gamma_c(D;A) \rightarrow \Gamma_c(X \times Y;A) \overset{\simeq}{\rightarrow} \Gamma_c(X;A Y) \]

For $p:D \rightarrow X \times Y $ be a separating object over $X \times Y$ as in the theorem. The isomorphism corresponds to the one of \ref{Prop_Gammac_exchange} when $X$ is $A$-separating, hence it is essentially the isomorphism of \ref{Lem_WeakseparatingTransport} applied to the the topos $\Tcal_{/X}$ and to the object $X \times Y$ by cosheaf extension from the separating objects of $\Tcal_{/X}$, but $D$ is one of the separating objects of $\Tcal_{/X}$ and hence one has a diagram (all the $\Gamma_c$ of the left square being computed in the topos $\Tcal_{/X}$):

\[\begin{tikzcd}[ampersand replacement = \&]
\Gamma_c(D;p_X^* A) \arrow{r} \arrow{d}{\simeq} \& \Gamma_c(X \times Y;p_X^* A) \arrow{d}{\simeq}  \arrow{r}{\simeq} \& \Gamma_c(X \times Y;A) \arrow{d}{\simeq}\\
\Gamma_c(X;p_X^*(A)D) \arrow{r} \& \Gamma_c(X, p_X^*(A)(X \times Y) ) \arrow{r}{\simeq} \& \Gamma_c(X;Y_A)
\end{tikzcd} \]

But the explicit description of the left most vertical arrow given in \ref{Lem_WeakseparatingTransport} allows to give an explicit description of the total diagonal map which is exactly the one presented in the theorem.

}

}

\block{
\Def{ We fix some set $(X_i)_{i \in I}$ of $A$-separating object of $\Tcal$, such that the $X_i$ and their sub-object form a generating family of $\Tcal$. And we define $\Ccal_c(\Tcal;A)$ to be the additive pseudo-category whose objects are the separating objects $X_i$ of $\Tcal$ and whose morphisms are the compactly supported map between the $A X_i$.
}

If the ground topos is locally decidable, one can find such a family formed of a single object $X$, in which case $\Ccal_c(\Tcal;A)$ will simply be a (non-unital) algebra. Because of the result above, this algebra should be thought of as the algebra of (finite) matrix with coefficients in $X$ in the sense that its elements corresponds to compactly supported function on $X \times X$.
But if we want to treat the case of a general basis we need a category.
}

\block{\Prop{Assume that $*$ is a linear involution on $A$ such that $(xy)^*=y^*x^*$ for all $x,y \in A$. Then $\Ccal(\Tcal;A)$ is a $*$-category for the $*$-operation which takes an arrow $f: X_A \rightarrow Y_A$ represented by a compactly supported function $f \in \Gamma_c(X \times Y,A)$ and exchange the variables of $f$ and apply $*$.
}

\Dem{This operation is clearly linear, we just have to check that $(fg)^*=g^* f^*$ but this follow easily from the description of a function $X_A \rightarrow Y_A$ in terms of a compactly supported section on $X \times Y$ given in proposition \ref{Prop_MatrixRep} using an easy computation very similar to the proof that $*$-transpose is a $*$-operation on matrix algebras with coefficient in a $*$-algebra.}

}

{\Prop{If $V$ is any sheaf of $A$-modules on $\Tcal$ then $X \mapsto \Kcal(A X, V)$ defines a non-degenerate right $\Ccal_c(\Tcal,A)$-module. }

\Dem{Let $X$ be one of the chosen $A$-separating objects. Let $f :X \rightarrow V$ be a compactly supported function, then by lemma \ref{Lem_cptSup_extension} one can find a compactly supported section $\lambda$ of $A$ on $X$ such that $\lambda$ is equal to $1$ on the support of $f$. Multiplication by $\lambda$ defines a compactly supported endomorphism of $X_A$ hence it is an element of $\Ccal_c(\Tcal,A)$ and $f \circ  \lambda = f$.}

}

\block{\Th{The functor $V \mapsto \Kcal(A X_i,V)$ defines an equivalence of categories between the category of sheaves of $A$-modules and the category of right non-degenerate $\Ccal_c(\Tcal,A)$-modules.}

\Dem{Let $A$-Mod denotes the category of sheaves of $A$-modules, let $\Ccal_c(\Tcal,A)$-Mod be the category of non-degenerate right module over $\Ccal_c(\Tcal,A)$, and  Let $S: A$-Mod $ \rightarrow \Ccal_c(\Tcal,A)$-Mod be the functor defined in the theorem. We will first construct a left adjoint $T$.

Let $M \in \Ccal_c(\Tcal,A)$-Mod, It suffices to construct an object $T(M)$ which satisfies the universal properties (naturally in $V$):

\[ \hom(T(M),V) = \hom(M,S(V)) \]

The functoriality of $T$ and the fact that it is adjoint of $S$ then follows from general categorical non-sense. A morphism $v$ from $M$ to $S(V)$ is the data for each $ X \in \Ccal_c(\Tcal,A)$ and each $m \in M(X)$ of a compactly supported morphism $v_{X,m}$ from $A X$ to $V$, which satisfies some equations translating the naturality of $v$. The key point is that the fact that the $v_{X,m}$ are compactly supported can be deduced from those equations. Indeed, as $M$ is non-degenerate, there exists an arrow $f \in \Ccal_{c} (\Tcal,A)$ such that $m=b.n$, hence by naturality $v_{X,m} = v_{Y,n} \circ f$ hence as $f$ is compactly supported $v_{X,m}$ is automatically compactly supported. Once this condition is removed a morphism in $\hom(M,S(V))$ is described as the data of functions from $A X$ to $V$ satisfying some relations, which can be translated as a map from some colimit $C$ to $V$, hence as the category $A$-mod is co-complete there is indeed such an object $T(M)$, which concludes the proof of the existence of the adjoint.

Now as $S$ commutes to co-limit (it follows immediately from proposition \ref{Prop_gammac_colimit}) $S \circ T(M)$ is defined by the same co-limit as $T(M)$ but computed in $\Ccal_c(\Tcal,A)$-Mod, hence it is isomorphic to $M$, moreover the unit of the adjunction $M \rightarrow S \circ T(M)$ is clearly an inverse of this isomorphism.

It remains to prove that the co-unit $c_N : T \circ S (N) \rightarrow N$ is also an isomorphism, but $S(C_N)$ is isomorphic to the identity of $S(N)$ hence it its enough to check that $S$ detect isomorphisms.

Let $f : N \rightarrow M$ be a map between two sheaf of $A$-modules such that $S(f)$ is an isomorphism. 
We will first prove that $f$ is monomorphism: let $v_1,v_2$ be two functions from  $V \subset X \in \Ccal_c(\Tcal,A)$ to $N$ such that $f \circ v_1 =f \circ v_2$. Let $U \ll V$ and let $\chi$ be a compactly supported function from $X$ to $A$ which is equal to one on $U$ and supported in $V$. The function $v_1 \chi$ and $v_2 \chi$ are compactly supported function from $X$ to $N$, hence the action of $f$ on them is the same as $S(f)$, hence $v_1 \chi = v_2 \chi$ hence $v_1 = v_2$ at least on $U$, but as this is true for any $U \ll V$ this show that $v_1= v_2$ on $V$ and as subobject of object in $\Ccal_{\Tcal,A}$ form a generating set by assumption it proves that $f$ is a monomorphism. A completely similar argument show that $f$ is also an epimorphism and this concludes the proof.
}
}

\blockn{When $A$ is the ring $\R$ or $\mathbb{C}$ of real or complex Dedekind numbers on $\Tcal$, the algebra $\Ccal_c(\Tcal,A)$ corresponds roughly to the convolution algebra of compactly supported function on a groupoid (see for example \ref{Ex_etalegpd}), as for groupoid algebras their is several way to complete it into a Banach algebra or $C^*$-algebra using different norm on this algebra. We will conclude this section by presenting the most important of these norms. For simplicity we will focus on the case where one use a single separating bound $X$ and hence that $\Ccal_c(\Tcal,\mathbb{R})$ is the algebra of endomorphisms of $X_{\mathbb{R}}$ with compact support on $X$. All the norm defined below (both internal and external) takes values en the upper semi-continuous real numbers (i.e. upper Dedekind cut).}

\block{We start with the $L^1$-norm or $I$-norm. Internally in $\Tcal$ the object $X_\mathbb{R}$ can be endowed with the $l^1$ norm:

\[ \Vert x \Vert_1 = \inf_{x=\sum \lambda_i x_i} |\lambda_i| \]

it is easy to check that this defines internally a pre-norm on $X_A$ and one can use it to put a operator norm on $\Ccal_c(\Tcal,\R)$, if $f \in \Ccal_c(\Tcal,\R)$:

\[ \Vert f \Vert_{I,l} = \sup_{x \in X_A, \Vert x Vert_1 \leqslant 1} \Vert f(x) \Vert_1 \]

in the sense that $\Vert f \Vert_{I,l}<q$ if and only if there is a $q'<q$ such that internally in $\Tcal$ one has $\forall x \in X_A$ such that $\Vert x \Vert_1 \leqslant 1$ one has $\Vert f(x) \Vert_1 \leqslant q'$.

One easily see that this is a norm on $\Ccal_c(\Tcal,\R)$ which satisfies $\Vert x y \Vert_{I,l} \leqslant \Vert x \Vert_{I,l}  \Vert x \Vert_{I,l} $. The involution is not isometric for this norm. To fix that one generally defines:

\[ \Vert f \Vert_I = \text{max} (\Vert f \Vert_{I,l},\Vert f^* \Vert_{I,l}) \]

Which is the so called $I$-norm or $L^1$ norm. The completion of $\Ccal_c(\Tcal,\R)$ for this norm is denoted $L^1(\Tcal,\R)$ and is a Banach algebra.

}

\block{One can try to define the $L^2$-norm or reduced norm in a similar way, but we need to deal with an additional difficulty: If one try to define the $l^2$ norm on $X_{\R}$ using the formula:

\[ \Vert x \Vert_2 = \left( \inf_{x=\sum \lambda_i x_i} |\lambda_i|^2 \right)^\frac{1}{2} \]

then this gives norm $0$ to all element of $X_A$, indeed, for any generator $x \in X$ one can write the corresponding element of $X_A$ as:
\[ x = \sum_{i=1}^n \frac{1}{n} x \]
and deduces that with the definition above the norm of $x$ is small than $\frac{1}{\sqrt(n)}$ and hence that it is zero. To avoid this, one need to add in the infimum defining $\Vert \Vert_2$ that the generators used in the expression are pairwise distinct but this work well only if $X$ is decidable. So one can construct the $L^2$-norm in a way similar to the $L^1$-norm above only if one can make our algebra to act on $X_A$ for $X$ a decidable object.

\bigskip

We proceed as follow: one choose $s:\Bcal \rightarrow \Tcal$ a surjection from a topos $\Bcal$ such that $s^*(X)$ is decidable in $\Bcal$. $\Bcal$ can for example be a boolean cover or the topos freely generated by adding a co-diagonal to $X$. Any $f \in \Ccal_c(\Tcal,\R)$ is then an endomorphism of $s^*(X)_{\R}$ and as $s^*(X)$ is decidable on can use this to define its $L^2$ operator norm.

Because the $l^2$ norm on $s^*(X)_{\R}$ is a Hilbert norm, the $L^2$-operator norm is preserved by the involution and satisfies the $C^*$-identity and the $C^*$-inequality $\Vert x^* x \Vert = \Vert x \Vert^2$ and $\Vert x \Vert ^2 \leqslant \Vert  x^* x +y^* y \Vert$. Hence the completion of $\Ccal_c(\Tcal,\R)$ for this norm is a real $C^*$-algebra called the reduced $C^*$-algebra of the topos and denoted $\Ccal^*_{red}(\Tcal,\R)$.
}

\block{Finally, one can define the maximal $C^*$-algebra of the topos as the universal real $C^*$-algebra $\Ccal^*_{max}(\Tcal,\R)$ generated by the involutive algebra $\Ccal_c(\Tcal,\R)$. In order to show that it exists one exactly need to prove the following lemma:

\Lem{For any $f \in \Ccal_c(\Tcal,\R)$ there exists a constant $K$ such that for any involutive morphism $h: \Ccal_c(\Tcal,\R) \rightarrow A$ into a $C^*$-algebra one has $\Vert h(f) \Vert \leqslant K$.}

Indeed, the max norm of $f$ is then defined as the infimum of all such constant $K$.

\Dem{We start with the case where $f$ is multiplication by a compactly supported function (also denoted $f$) on $X$, and assume that $O \leqslant f \leqslant 1$. 

Any positive function on $X$ can be written as $g^* g$ so for any morphism $h$ it is sent to a positive element of the $C^*$-algebra $A$. Hence in this case $h(f)$ is a positive element, for the same reason if $g \leqslant f$ are two positive functions on $X$ then $h(g) \leqslant h(f)$. In this case, as $f \leqslant 1$ one has $f^2 \leqslant f $ and hence $h(f)^2 = h(f^2) \leqslant h(f)$ which proves that $\Vert h(f) \Vert \leqslant 1$, hence $K=1$ works.

For a general compactly supported function $f$ on $X$, if $K$ is a constant larger than $|f|$ then $f^*f/K^2$ is a positive function between $0$ and $1$ hence for all morphism $h$ one has $\Vert h(f) \Vert \leqslant K$. So this $K$ works for $f$.

Let now $U \subset X$ be any sub-object, $s: U \rightarrow X$ be any map, and $\lambda$ be a compactly supported functions in $U$. One can then define an element of $s \lambda \in \Ccal_c(\Tcal,\R)$ as follow, internally as a function from $X$ to $X_A$ by:

\[ s\lambda(x)= \lambda(x) s(x) \]

Which is compactly supported (its support is the support of $h$ which is compact in $U \subset X$) One easily check that his adjoint is:

\[ (s\lambda)^* (x) = \sum_{s(y)=x} \lambda(y) y \]

and that $(s\lambda)(s\lambda)^*$ satisfies:

\[ (s\lambda)(s\lambda)^*(x) = \left(\sum_{s(y)=x} \lambda(y)^2 \right) x \] 

hence $(s\lambda)(s\lambda)^*$ is multiplciation by a (compactly supported) function, for any morphism $h$ one has indeed a constant $K$ such that $\Vert h((s\lambda)(s\lambda)^*) \Vert \leqslant K$ and hence $\Vert h(s \lambda) \Vert \leqslant K$, and one also have a $K$ for this type of elements.

Take now a general elements $f$ of $\Ccal_c(\Tcal,\R)$, i.e. a compactly supported function from $X_A$ to $X_A$.

Using theorem \ref{Prop_MatrixRep}, the morphism $f$ can be represented by an $\lambda_0$ element of $\Gamma_c(X \times X,A)$, moreover, as $X$ is a bound there is a finite co-product of subobject $U_i$ of $X$:

 \[ D=\coprod_{i=1}^n U_i \]

endowed with a map $(p_1,p_2) : D \rightarrow X \times X$ and a compactly supported function $\lambda \in \Gamma_c(D,A)$ such that $\lambda_0$ is the image of $\lambda$ by the map $\sigma_{(p_1,p_2)}$, following theorem \ref{Prop_MatrixRep}, this mean that $f$ can be described as:

\[ f(x) = \sum_{p_1(d)=x} \lambda(x) p_2(x) \]

If $D$ is just one sub-object $U \subset X$, then the above formula for $f$ can be re-written as:

\[ f =  (p_2 \lambda_2) (p_1 \lambda_1)^* \]

where $\lambda_2$ and $\lambda_1$ are two compactly supported functions on $U$ such that $\lambda= \lambda_1 \lambda_2$.
In the more general case where $D$ is indeed a co-product of $U_i \subset X$ then it corresponds to a decomposition of $f$ as:

\[ f = \sum_{i=1}^n (p^i_2 \lambda^i_2) (p^i_1 \lambda^i_1)^* \]

but we already know that one has constants $K^i_1$ and $K^i_2$ that control the norm of $h((p^i_1 \lambda^i_1)^*)$ and $h(p^i_2 \lambda^i_2)$ for any morphism $h$ as above and hence one has:

\[ \Vert h(f) \Vert \leqslant \sum K^i_1 K^i_2 \]

which concludes the proof.

}

}

\section{Examples}
\label{Sec_examples}

\blockn{We will conclude this paper by giving some example of topos to which the above formalism apply and the corresponding algebras. Some of the following examples are only sketched.}

\block{Let $\Lcal$ be a Hausdroff (equivalently regular) locally compact locale. Then $\Lcal$ is absolutely locally compact as a topos, the terminal object is separating and the subobject of the terminal object form a basis of quasi-finite object of cardinal $1$.

Hence any c-soft sheaf of rings on $X$ is admissible. Moreover the family $X_i$ can be chosen reduced to the single object $1$. It can also be check that conversely any admissible sheaf of ring have to be c-soft.

hence our main result become:

\Prop{Let $X$ be a locally compact Hausdorff locale and $A$ a c-soft sheaf of ring over $X$. The category of sheaves of $A$-modules is equivalent to the category of non-degenerate modules over the ring $\Gamma_c(A)$ of compactly supported sections of $A$.}

The theorem applies in particular when $A$ is the sheaf of continuous function with value in $\mathbb{R}$ or $\mathbb{C}$ (or any unital topological $\mathbb{R}$-algebra), as soon as we assume either the axiom of dependant choice or that $X$ is completely regular. One can of course take $X$ to be any locally compact Hausdorff topological space or any Hausdorff manifold (in which case one does not need the axiom of choice).

If $A$ is endowed with an involution (for example the identity if $A$ is commutative or complex conjugation if $A = \mathbb{C}$ then the involution on $\Gamma_c(A)$ is just the ``pointwise'' application of the involution on $A$.

}

\block{Let now $\Lcal$ be a locale which is only locally a locally compact Haussdorff locale.
Let $(U_i)_{i \in I}$ be an open covering of $\Lcal$ by open sublocales which are Hausdorff and locally compact. For each $i$, $\Lcal_{/U_i}$ is just $U_i$ and is a locally compact Hausdorff locale. In particular $\Lcal$ is absolutely locally compact and each $U_i$ as well as $\coprod_{i \in I} U_i$ (if $I$ is decidable) form separating objects.

As in the previous example any sheaf of rings $A$ which is c-soft on each $U_i$ will be admissible with the $U_i$ $A$-separating and one can take the $(U_i)_{i \in I}$ as our family $X_i$.

\Prop{In this situation, with $A$ a sheaf of ring which is c-soft of each $U_i$,

\begin{itemize}
\item  $\Ccal_c(\Lcal,A)$ can be chosen to be the additive pseudo-category whose objects are the $i \in I$ and whose morphism from $i$ to $j$ are compactly supported sections of $A$ on $U_i \wedge U_j$ composition is just the multiplication in $A$.

\item If $A$ is involutive the involution is just pointwise application of the involution of $A$ and exchange of the source and the target.

\item if $I$ is decidable (or if we assume the law of excluded middle), one can take $\coprod U_i$ as the only object in the familly $X_i$. In this case $\Ccal_c(\Lcal, A)$ is the algebra of finitely supported matrix with coefficient in $I$ whose $i,j$ component is a compactly supported section of $A$ on $U_i \wedge U_j$. Multiplication being define by matrix multiplication and the multiplication in $A$. 

\item In this case the involution is matrix transposition together pointwise application of the involution in $A$.

\item In both case, non-degenerate modules over $\Ccal_c(\Lcal,A)$ are the same sheaf of $A$-module over $\Lcal$.

\end{itemize}
}

Here again this applies when $A$ is the sheaf of real or complexe valued continuous functions as soon as we assume the axiom of dependant choice or that each $U_i$ is completely regular.

}

\block{\label{Ex_etalegpd}Let now $G=(G_0,G_1,s,t,\mu)$ be an étale topological groupoid. Let $\Tcal$ be the topos of $G$-equivariant sheaves, i.e. sheaves over $G_0$ endowed with an action of $G$ and $G$-equivariant maps between them.

The forgetful functor from $G$-equivariant sheaves to sheaves over $G_0$ is the $f^*$ part of an étale surjective geometric morphism $f: G_0 \rightarrow \Tcal$ which corresponds to the object $X$ which is $G_1$ over $G_0$ endowed with its multiplication action on itself.

Moreover $G_1$ can be described as $G_0 \times_{\Tcal} G_0$ and the groupoid structure on $G$ can be recovered as the obvious ``pair groupoid'' structure coming from this description of $G_1$. Any topos admitting an étale covering by a locale $G_0$ is of this form, such topos is called an étendu.

If $G_0$ is locally compact and Hausdorff, then any ring object $A$ of $\Tcal$ which is c-soft when seen as a sheaf over $G_0$ is admissible with $X$ as a $A$-separating object. The sub-object of $X$ form a generating family of $\Tcal$, hence one can take $X$ as the single element of the family $(X_i)_{i \in I}$.

By theorem \ref{Prop_MatrixRep}, compactly supported endomorphisms of $X_A$ are the same as compactly supported sections of $A$ on $X \times X$ but $\Tcal_{X \times X}$ is exactly the locale $G_1$ hence we need to distinguishes two cases:

Either $G_1$ is Hausdorff in this case ``compactly supported sections'' of $X$ do mean sections which are compactly supported. Or $G_1$ is not Hausdorff, in which case compactly supported sections on $G_1$ are computed using the co-sheaf construction of \ref{Th_def_cptsection}. As $G_1$ is always locally Hausdorff (it is étale over $G_0$ which is Hausdorff), compactly supported sections on $G_1$ can be computed using the cosheaf property on the covering of $G_1$ by some Hausdorff open sub-objects, and will be exactly linear combinations of compactly supported functions on Hausdorff open subspaces as it is usual in non-comutative geometry. In both case it is not very hard to see that multiplication and involution are the usual multiplication and involution of étale groupoid algebra.
}

\block{In the special case of an étale groupoid whose space of object is locally compact Hausdorff and totally disconnected (i.e. with a basis of compact open subspaces stable under intersection) then any sheaf over the space of objects is c-soft because any open can be covered by compact clopen subspaces. Hence in this case any sheaf of ring is admissible. Applying the above machinery to such a groupoid with a constant sheaf of rings gives exactly the algebra constructed by B.Steinberg in \cite{steinberg2010groupoid} and the equivalence between modules on the algebra and sheaf of module over the topos is the main result of his second paper on the subject \cite{steinberg2014modules}.}

\block{Let now look at a simple example where the divisibility condition in the definition of admissible sheaf of ring is not vacuous. Let $G$ be a pro-finite group, in fact, for simplicity, take $G = \Z_p$ the additive group of p-adic integer.
Let $\Tcal$ be the topos of smooth $G$-set, i.e. the category of sets endowed with a continuous (i.e. smooth, or locally constant) action of $\mathbb{Z}_p$. This is a absolutely compact topos and $1$ is separating. A ring object in $\Tcal$ will be admissible if and only if $p$ is invertible in $A$. Indeed, as the topos is atomic the softness condition is vacuous (all the $Loc(X)$ are discrete, so one always have compactly supported functions) but the quasi-finite generators corresponds to $G$-orbits whose cardinal is always a power of $p$, so one need to have an inverse for $p$. 

Again, for simplicity we will focus on the case where $A$ is a ring with trivial $\mathbb{Z}_p$ action and in which $p$ is invertible. We take the $(X_k = \mathbb{Z}/p^k \Z)$ as our family of generators. A function from $X_k$ to $X_{k'}$ in $\Ccal(\Tcal,A)$ is a compactly supported section on $X_k \times X_{k'}$ the orbits of $X_k \times X_{k'}$ corresponds to the double cosets $(p^{k'}\Z_{p}) \setminus \Z_p /(p^k \Z_{p})$, so morphisms from $X_k$ to $X_{k'}$ corresponds to linear combinations of this with coefficient in $A$ and the composition can be seen to be the multiplication of double cosets algebra.

Note that in this case the invertibility of $p$ appears to be unimportant for the definition of the algebra but are important for the proof that modules over this algebra are the same as module objects in $\Tcal$. The main result relating module objects over $A$ and modules over $\Ccal_c(\Tcal,A)$ is, in this case, the well known relation between representations of the double cosets algebra and representations of the group.
}

\block{\label{example_Graph}We now consider a finite\footnote{having for each verticies a finite number of edges targeting it would be enough, but it is simpler to assume the graph is finite for certain details below.} directed graph $\G = (\G_0,\G_1)$, i.e. $\G_0$ is a finite set (its elements are called vertices), $\G_1$ is a finite set whose objects are called arrows, and there is two maps $s,t : \G_1 \rightrightarrows \G_0$ giving respectively the source and the target of each arrow.

\bigskip

If $\G$ is a graph a $\G$-presheaf $\Fcal$ is the data of:

\begin{itemize}

\item For each vertices $x \in \G_0$ a set $\Fcal(x)$.

\item For each arrow $a \in \G_1$, $a:x \rightarrow y$, i.e. $x=s(a)$ and $y=t(a)$, one has a function $\Fcal(a):\Fcal(y) \rightarrow \Fcal(x)$.

\end{itemize}
Morphisms of $\G$-presheaves are naturally defined as collection of functions $f_x: \Fcal(x) \rightarrow \Fcal'(x)$ for $x \in \G_0$ such that for all $a \in \G_1$ the induced square commute.

A $\G$-sheaf is a $\G$-presheaf such that for all $x \in \G_0$ one has:

\[ \Fcal(x) \simeq \prod_{a: y \rightarrow x \in \G_1} \Fcal(y) \]

Where the isomorphism is induced by the natural map which is $\Fcal(a) : \Fcal(x) \rightarrow \Fcal(y)$ on the component $a :Y \rightarrow x$.

The category of $\G$-presheaf is equivalent to the category of presheaf on the category $\G^p$ freely generated by $\G$ i.e. the category of paths in $\G$. We will show that the $\G$-sheaves are the sheaves for a Grothendieck topology on $\G^p$.

For any verticies $x$ of $\G$, let $x^-$ be the sieve over $x$ of morphism (i.e. path) from $y$ to $x$ which are of length at least $1$. Equivalently, $x^-$ is the sieve generated by the covering family of the $y \rightarrow x$ for all arrows to $x$ in the graph. if $y \rightarrow x$ and $y' \rightarrow x$ are two arrow in the graph it is not hard to see that the pullback of these two arrows (in the category of presheaf) is $\emptyset$ unless they are equal in which case it is the arrow itself. In particular the sheaf condition with respect to the sieve $x^- \hookrightarrow x$ is exactly the condition that the map:

\[ \Fcal(x) \simeq \prod_{a: y \rightarrow x \in \G_1} \Fcal(y) \]

is an isomorphism.

The pullback of $x^{-} \hookrightarrow x$ by any morphism $y \rightarrow x$ is the maximal sieve $y \hookrightarrow y$ as soon as the morphism $Y \rightarrow x$ has length at least one, and is $y^{-} \hookrightarrow y$ if it is the identity. In particular, the family of all sieve $x^{-} \hookrightarrow x$ and $x \hookrightarrow x$ is stable under pullback. In particular, by \cite[Corrolary II.2.3]{SGA4I}, in order to check whether a pre-sheaf is a sheaf with respect to the topology generated by the covering sieve $x^{-} \hookrightarrow x$ one just have to check to sheaf condition for those sieves, i.e. $\G$-sheaves are exactly the sheaves for the topology generated by the $x^{-} \hookrightarrow x$. In particular, $\G$-sheaves form a topos.

\bigskip

We will denote by $\Tcal_{\G}$ the topos of $\G$-sheaves. For example, if $G$ has one verticies and one arrow it is exactly the topos $B\mathbb{Z}$ of sets with an action of $\mathbb{Z}$, if $G$ has one verticies and $2$ arrows it is the so-called ``J\'onsson-Tarski'' topos of sets $X$ endowed with an isomorphism between $X$ and $X \times X$. The reader should immediately notes that the site we used to define it is far from being sub-canonical: if one fix a verticies $y$ then $\Fcal(x):=\{\text{Path from $x$ to $y$ }\}$ does not satisfies the sheaf condition.

\bigskip

Let $x$ be a verticies of the graph. Let $P_x$ the representable sheaf associated to the object of $\G^p$, i.e. the sheafification of the representable pre-sheaf. The topos $\Tcal_{\G}/P_x$ can be described by the slice site  $\G^p/x$. As a category, it is the poset of finite paths in $\G$ that ends at $x$ ordered by extension, the topology is generated by the cover $p^-$ of a path $p$ by all the paths of length one more than the length of $p$ which extend $p$. This is exactly a site for the locale of infinite path in $G$ ending at $x$ (i.e. indexed by $i \leqslant 0$ and such that $p_0=x$ ).

In particular, $\Tcal_{\G}/P_x$ is a locale and even a Stone space. Hence, $\Tcal_{G}$ is absolutely locally compact, $P_X$ is a separating object for $\Tcal_{\G}$, any sheaf of ring over $\Tcal_{\G}$ is admissible, and compactly supported section over $P_x$ are just ordinary section.

\bigskip

Let $K$ be a unital ring, and consider it as a constant sheaf over $\Tcal_{\G}$. One take $X= \coprod_{x \in \G_0} P_x$ as a single generator of the topos to construct the algebra $K_{\G} = \Ccal_c(\Tcal_{\G},K)$.

If $M$ is a $K$ module over $\Tcal_{G}$ then it is just a $\G$-sheaf of ordinary $K$-modules, and the corresponding $\Ccal_c(\Tcal_{\G},K)$-module is (at least at the level of the underlying $K$-module) :

\[ \bigoplus_{x \in \G_0} M(x) \]

\bigskip

\Def{The Leavitt path algebra $L_K(\G)$ is the involutive $K$-algebra generated by elements $v_a$ for $a \in \G_1$ and $p_v$ for $v \in \G_0$ with the relation:

\begin{itemize}

\item $p_e^*=p_e$ and $p_e p_{e'} = \delta_{e,e'} p_e$.

\item $v_a p_{s(a)} =p_{t(a)} v_a = v_a$

\item $v_a^* v_b = \delta_{a,b} p_{s(a)}$

\item $ p_e = \sum_{x \in \G_1, t(x)=e} v_x v_x^* $

\end{itemize}

}

One will first see that a (non-degenerate) right $L_K(\G)$-module is exactly a $K$-module over $\Tcal_{\G}$: First observe that $i=\sum_{e \G_0} p_e$ is a unit for this algebra. Indeed it follows easily from the relation above that $ip_e=p_e i = p_e$ and $i v_a = v_a i =v_a$ and $iv_a^* = v_a^* i = v_a^*$ and a general element is a polynomial in those elements hence $i$ is a unit. In particular the $p_e$ form a maximal family of orthogonal projections hence correspond to decomposition of $M$ into:

\[ M = \bigoplus M(e) \]

where $M(e) = M.p_e$.

Moreover, right multiplication by $v_a$ for $a:e \rightarrow e'$ an arrow in $\G$ corresponds exactly to a linear map from from $M(e')$ to $M(e)$, indeed for an element $x$ in $M(e'')$ with $e'' \neq e'$, one has $x v_a= x p_{e''} v_a = 0$ and for $x \in M(e')$ one has $x v_a p_e = x v_a$ hence $x v_a \in M(e)$. If one considers the $K$-algebra generated by the $v_a$ and $p_e$ subject to the first two relations, a right module over this algebra would correspond exactly to a $\G$-presheaf of $K$-module. Adding the existence of the $v_a^*$ subject to the last two relations exactly assert (in terms of by-product) that the natural map:

\[ M(e) \rightarrow \prod_{ a: e' \rightarrow e} M(e') \]

is an isomorphism and hence a right $L_K(\G)$-module can be identified with a $\G$-sheaf of $K$-module

\bigskip

A map from $\coprod_{e \in \G_0} P_e$ to a sheaf of $K$-modules $M$ corresponds exactly to an element of $\prod M(e)$ hence to an element of the $L_K(\G)$-module corresponding to $M$, hence the free $K$-module on $\coprod_{e \in \G_0} P_e$ corresponds to free $L_K(\G)$-module on one generator. The algebra of the topos $\Tcal_{\G}$ is the algebra of compactly supported endomorphisms of this free module, but because of the compactness of $\coprod_{e \in \G_0} P_e$ it is exactly the set of all endomorphisms of $L_K(\G)$ hence it is $L_K(\G)$ itself. One has proved that:

\Th{The convolution algebra $\Ccal_c(\Tcal_{\G},K)$ is the Leavitt path algebra $L_K(\G)$.}

The result still hold in the case where the graph $\G$ is infnite as long as any vertices has only a finite number of edges pointing to it. In this infinite situation, the algebra is no longer unital but the finite sum of $p_e$ form an ``approximate unit'' and not all endomorphisms of $L_K(\G)$ are compactly supported but it is not very hard to make everything works. Also in a constructive context it is useful to assume that the set $\G_0$ of verticies of $\G$ is decidable otherwise one cannot consider the object $\coprod_{e \in \G_0} P_e$ as a separating object and one need to work with a ``Leavitt path algebroid'' whose set of object is $\G_0$ instead.

}

\bibliography{Biblio}{}
\bibliographystyle{plain}

\end{document}